\documentclass[12pt,a4paper,twoside]{article}
\usepackage[a4paper,margin=2.5cm,bottom=4cm]{geometry}			

\RequirePackage{fix-cm}

\usepackage{natbib}
\usepackage{todonotes}
\usepackage{graphicx}
\usepackage{amsfonts}
\usepackage{hyperref}
\usepackage{amsmath} 
\usepackage{amssymb} 
\usepackage{mathabx} 
\usepackage{mathtools}

\usepackage{pdflscape}

\usepackage{caption}
\usepackage{subcaption}

\usepackage{booktabs}

\usepackage{multirow}
\usepackage{breakcites}

\begin{document}

\title{
	Obey validity limits of data-driven models
	\thanks{This work was supported by the Deutsche Forschungsgemeinschaft (DFG, German Research Foundation) under Germany's Excellence Strategy - Cluster of Excellence 2186 ``The Fuel Science Center''.}
}

\title{Obey validity limits of data-driven models}   

\author{Artur M. Schweidtmann$^{1,*}$ \and Jana M. Weber$^2$ \and Christian Wende$^1$ \and Linus Netze$^1$ \and Alexander Mitsos$^{3,1,4}$}

\maketitle \noindent
$^1$ Process Systems Engineering (AVT.SVT), RWTH Aachen University, Aachen, Germany. \\
	$^2$ Department of Chemical Engineering and Biotechnology, University of Cambridge, Cambridge, United Kingdom. \\
	$^3$ JARA-CSD, 52056 Aachen, Germany. \\
	$^4$ Institute of Energy and Climate Research, Energy Systems Engineering (IEK-10), Forschungszentrum J\"ulich GmbH, 52425 J\"ulich, Germany. \\ 
	$^*$ artur.schweidtmann@rwth-aachen.de \\[2ex]%

\sloppy

\begin{abstract}
Data-driven models are becoming increasingly popular in engineering, on their own or in combination with mechanistic models. 
Commonly, the trained models are subsequently used in model-based optimization of design and/or operation of processes. 
Thus, it is critical to ensure that data-driven models are not evaluated outside their validity domain during process optimization.
We propose a method to learn this validity domain and encode it as constraints in process optimization.
We first perform a topological data analysis using persistent homology identifying potential holes or separated clusters in the training data.
In case clusters or holes are identified, we train a one-class classifier, i.e., a one-class support vector machine, on the training data domain and encode it as constraints in the subsequent process optimization.
Otherwise, we construct the convex hull of the data and encode it as constraints.
We finally perform deterministic global process optimization with the data-driven models subject to their respective validity constraints. 
To ensure computational tractability, we develop a reduced-space formulation for trained one-class support vector machines and show that our formulation outperforms common full-space formulations by a factor of over 3,000, making it a viable tool for engineering applications.
The method is ready-to-use and available open-source as part of our MeLOn toolbox (\url{https://git.rwth-aachen.de/avt.svt/public/MeLOn}).
\end{abstract}

\section{Introduction}
\label{intro}
Supervised machine-learning techniques have been re-emerging as a promising avenue for data-driven modeling in various engineering disciplines~\citep{venkatasubramanian2019promise}.
In most applications, the overall goal is the optimal decision-making based on available data and \textit{a priori} knowledge.
Thus, data-driven models and mechanistic models are often combined to form hybrid models~\citep{mogk2002application,kahrs2008incremental,von2014hybrid,glassey2018hybrid}.
Subsequently, hybrid models are frequently used in model-based optimization of design and/or operation of processes~\citep{mcbride2019overview,Schweidtmann2019detGlobalANN}.  \newline \indent
A critical issue of data-driven models is their limited extrapolability.
Unless strong assumptions are posed on the learned function, data-driven models can only be valid in regions where they have sufficiently dense coverage of training data points.
We refer to this as the \textit{validity domain}~\citep{courrieu1994three}.
Consequently, there is a need to avoid the evaluation of data-driven models outside their validity domain during optimization. 
Note that we refer to the validity domain of individual data-driven models throughout this work, but the concept can also be applied to hybrid models~\citep{kahrs2007validity}. \newline \indent
The vast majority of previous publications use box constraints (i.e., hyperrectangles) to bound the inputs of data-driven models, i.e., each variable has independent bounds.
This approach is practical when the training data is obtained from simulations based on regular grids or Latin hypercubes that are sufficiently dense.
It is also advantageous for local and global optimization.
However, it requires \textit{a priori} known bounds and the possibility to obtain training data for any input combination. 
In practice, simulations can fail \citep{asprion2020modeling} and industrial process data usually does not cover hyperrectangular spaces~\citep{Asprion2019GrayBox}. 
This leads to manual selection of wrong bounds, which may cut off optimal solutions or overestimate the validity domain. \newline \indent
As proposed by \citet{courrieu1994three}, a few previous works in process systems engineering (PSE) constructed the convex hull of the training data points to describe the validity domain and integrated it as a set of linear constraints in optimization problems~\citep{kahrs2007validity,Zhang.2016b,Asprion2019GrayBox}. 
By definition, the convex hull is the smallest convex set that contains all data points.
Commonly, evaluations of data-driven models inside the convex hull of the training data are called interpolation and outside extrapolation. 
However, roughly speaking, the convex hull cannot distinguish between for potential holes in the training data set, gaps between separated clusters of data, and nonconvex boundaries. 
Thus, staying within the convex hull seems only a necessary condition for the validity of data-driven models and not sufficient. 
Identifying if the convex hull is a suitable model for the data domain is very challenging in high dimensions. 
Notably, \citet{Zhang.2016b} extended the convex hull method to the union of multiple polytopes by introducing binary variables and additional constraints to the problem. 
However, this algorithm becomes impractical when the number of data points or their dimension is very high.
Besides convex hull formulations, there exist also several publications that circumvent extrapolation problems en passant in different ways. 
For instance, \citet{mistry2018optimization} penalize deviations from a training data mean in a space that is parameterized using principal component analysis.
\citet{rall2019rational} constrain the maximal allowed distance from the nearest training data point resulting in a nonsmooth optimization problem.
\citet{kumar2019machine} train multiple data-driven models on a design problem and reject designs where the variation between the models is large.
Similarly, \citet{pinto2019bootstrap} use bootstrap aggregation to estimate error bounds for hybrid mechanistic/data-driven models.
There exist further methods that quantify the variance or confidence interval of predictions such as Bayesian methods and maximum likelihood estimations \citep{papadopoulos2001confidence}.
However, this leads to chance-constrained programming problems~\citep{charnes1959chance,schweidtmann2020globalGP}. 
Likewise, there are a few studies on the adaptive exploration of the design space~\citep{larson2012design,chen2018optimization,knudde2019active} and related works on constrained Bayesian optimization~\citep{Shahriari.2016}.
However, we focus on fixed training data sets in this study while the extension of our method to adaptive sampling is a promising future research. \newline \indent 
An alternative to box constrains and convex hull is to use a nonlinear classifier that can also model complicated validity domains.
A few previous studies in mechanical engineering~\citep{malak2010using,roach2011improved} used Support Vector Domain Description (SVDD) \citep{tax1999data} to model the validity domain of data-driven equipment models.
Also, \cite{quaglio2018model} use binary support vector classification to include reliability constraints into model-based design of experiment.
As only valid training data points are given in most engineering applications, we consider one-class classification in this work. 
There exists a broad variety of one-class classifiers that can be divided into density methods, boundary methods, and reconstruction methods~\citep{Tax2001}.
Also, one-class classification is closely related to novelty, outlier, or anomaly detection~\citep{chandola2009anomaly,pimentel2014review,khan2009survey,khan2014one,ding2014experimental}. 
The previous literature indicates that one-class support vector machines (SVMs)~\citep{scholkopf2000support} are common and suitable for the problem at hand. 
Compared to density models, less training data is required to construct the boundary, since only the boundary is estimated and not a complete density distribution~\citep{Tax2001}. 
In addition, the one-class SVM is tolerant to outliers in the training set~\citep{pimentel2014review}. \newline \indent
Optimization problems with one-class SVMs embedded are nonconvex.
Thus, deterministic global optimization is desirable to identify global solutions.
However, these models lead to large-scale optimization problems that are difficult to solve.
In our previous work, we showed that a reduced-space (RS) formulation and the use of McCormick relaxations are advantageous for the optimization of two other important classes of data-driven models, namely artificial neural networks~\citep{Schweidtmann2019detGlobalANN} and Gaussian processes~\citep{schweidtmann2020globalGP}.  
We propose a similar idea here for one-class SVM. \newline \indent
The global shape of data matters because it often provides important information about the underlying phenomena represented by the data. 
Especially in high-dimensional data, topological data analysis~(TDA) can reveal and quantify objects and features not directly visible to the human eye.  
In the context of this work, it provides valuable information about typologies in the training data that can be colloquially thought of as holes or separated clusters.
TDA was initiated relatively recently~\citep{letscher2002topological,zomorodian2005computing}.
Its roots lie in applied (algebraic) topology and computational geometry~\citep{chazal2017introduction} and it is commonly used to account for higher-order interactions in data, to comprehend mesoscale structures, or to compare different data spaces~\citep{patania2017topological}.
The most common TDA method is persistent homology~\citep{wasserman2018topological}. 
So far, there are only a few applications of persistent homology in the fields of (bio-)chemical engineering and material science~\citep{hiraoka2016hierarchical,saadatfar2017pore,xia2018persistent,xia2019persistent,smith2020topological}.
\newline \indent
We propose a three-step approach to model the validity domain of data-driven models for optimization.
We first perform TDA using persistent homology.
In case clusters or holes are identified, we train a one-class SVM on the training data domain of the data-driven models and encode it as constraints in the subsequent process optimization. 
Otherwise, we construct the convex hull of the data and encode it as constraints.
We finally perform deterministic global process optimization with the data-driven models and their respective validity constraints. 
To ensure computational tractability, we develop a RS formulation for trained one-class SVMs. 
Moreover, we employ convex and concave envelopes of kernel functions to accelerate optimization. 
We demonstrate the potential of our method on a set of illustrative mathematical case studies and an engineering case study, i.e., the open-loop control of a sulfur recovery unit. 
\section{Methodology}
\begin{figure}
	\centering
	\includegraphics[width = .45\textwidth]{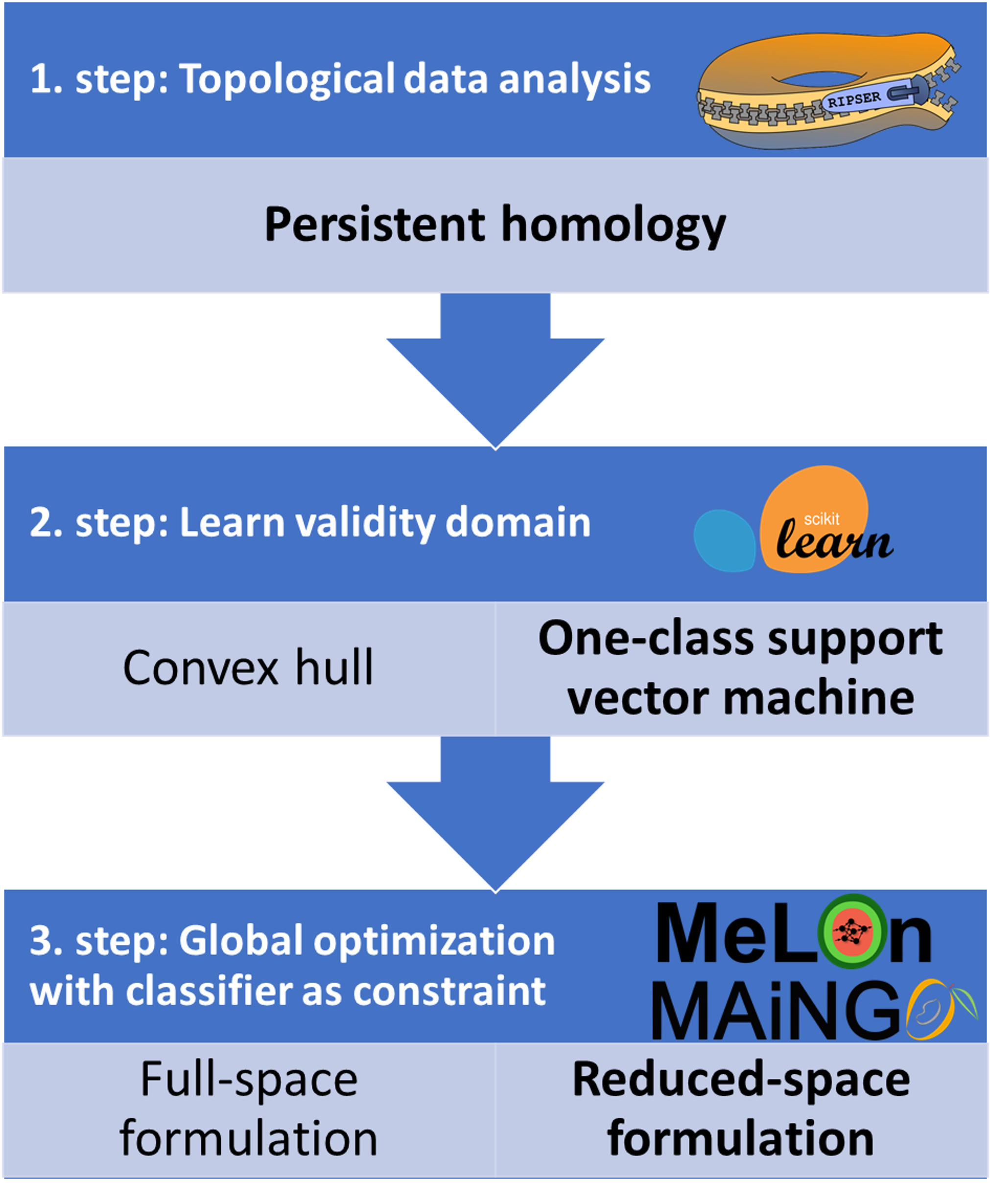}
	\caption{Overview of the proposed three step methodology to obey validity limits of data-driven models in optimization}
	\label{fig:Methodology_Overview}
\end{figure}

As illustrated in Figure~\ref{fig:Methodology_Overview}, we propose a three step approach to obey validity limits of data-driven models during optimization.
In the first step, we conduct a TDA of the training data.
In the second step, we either construct the convex hull of the data or we train a one-class classifier, i.e., a SVM.
In the third step, we embed the trained classifier or convex hull in the optimization problem and solve it to global optimality. 
The described methods are available open-source.
We use the \textit{Ripser.py} toolbox that is available open-source under MIT license in Python for performing the TDA~\citep{ctralie2018ripser}. 
The training of the one-class SVM is performed by Scikit-learn~\citep{pedregosa2011scikit} and the convex hulls are identified using SciPy~\citep{SciPy}.
We provide the one-class SVM  within the ``MeLOn - \textbf{M}achin\textbf{e} \textbf{L}earning Models for \textbf{O}ptimizatio\textbf{n}'' toolbox under the Eclipse public license~\citep{MeLOn_Git}.
The resulting optimization problems are solved using our open-source global solver MAiNGO~\citep{MAiNGO}.
\subsection{Topological data analysis using persistent homology}
\label{sec:Lear_validity_domain_Background_Topological_Data_Analysis}
\begin{figure*}
	\centering
	\includegraphics[width=\textwidth]{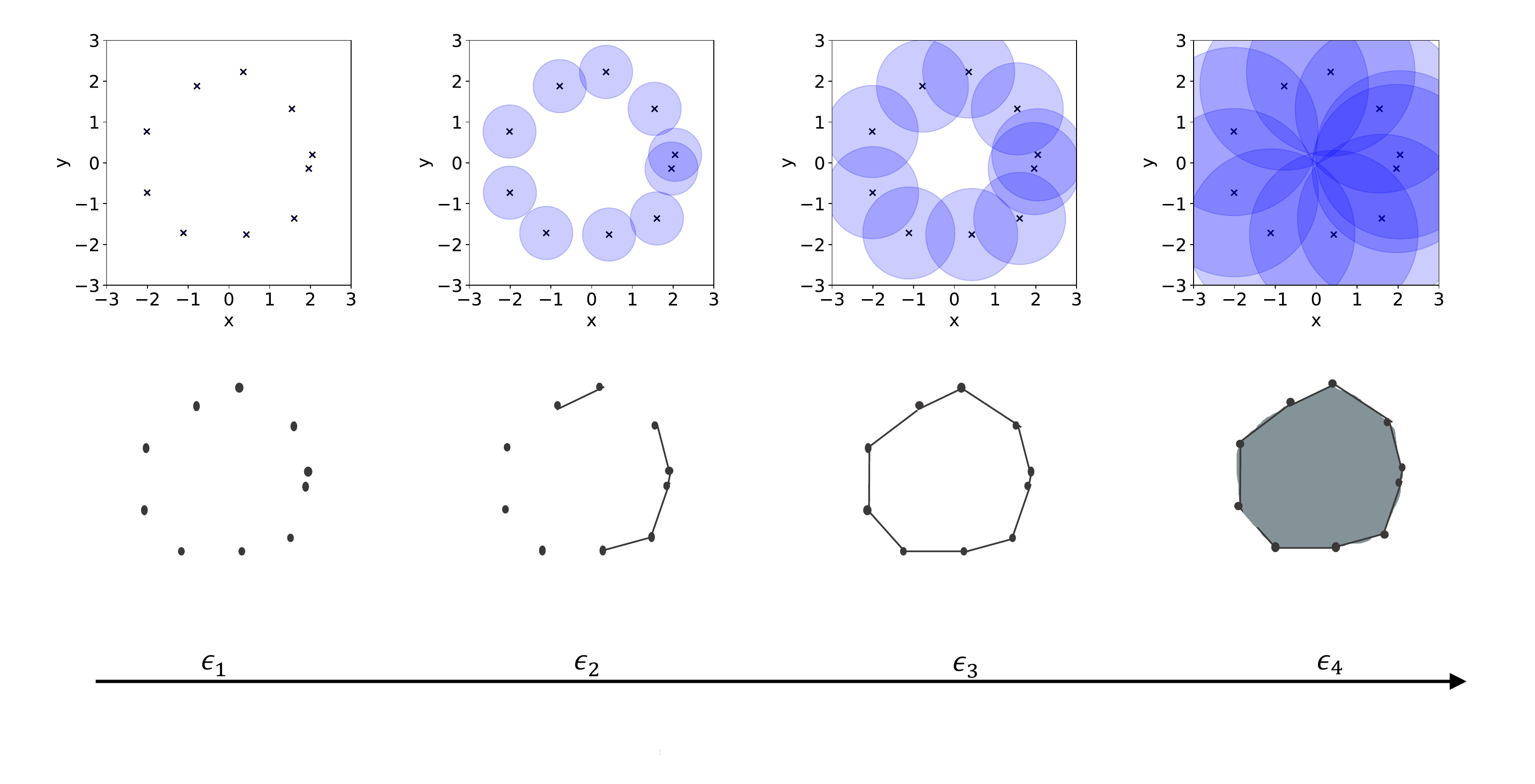}
	\caption{Illustration of a Vietoris-Rips filtration utilized for persistent homology. 
		The upper part shows the data set and circles around the data points with increasing diameter $\epsilon$. 
		The bottom image illustrates the simplicial complexes formed during the filtration. 
		The figure is based on \cite{kimura2017quantification}}
	\label{fig:LVD_Method_Illustration_Rips_filtration}
\end{figure*}

\begin{figure} 
	\centering 
	\includegraphics[width=0.5\textwidth]{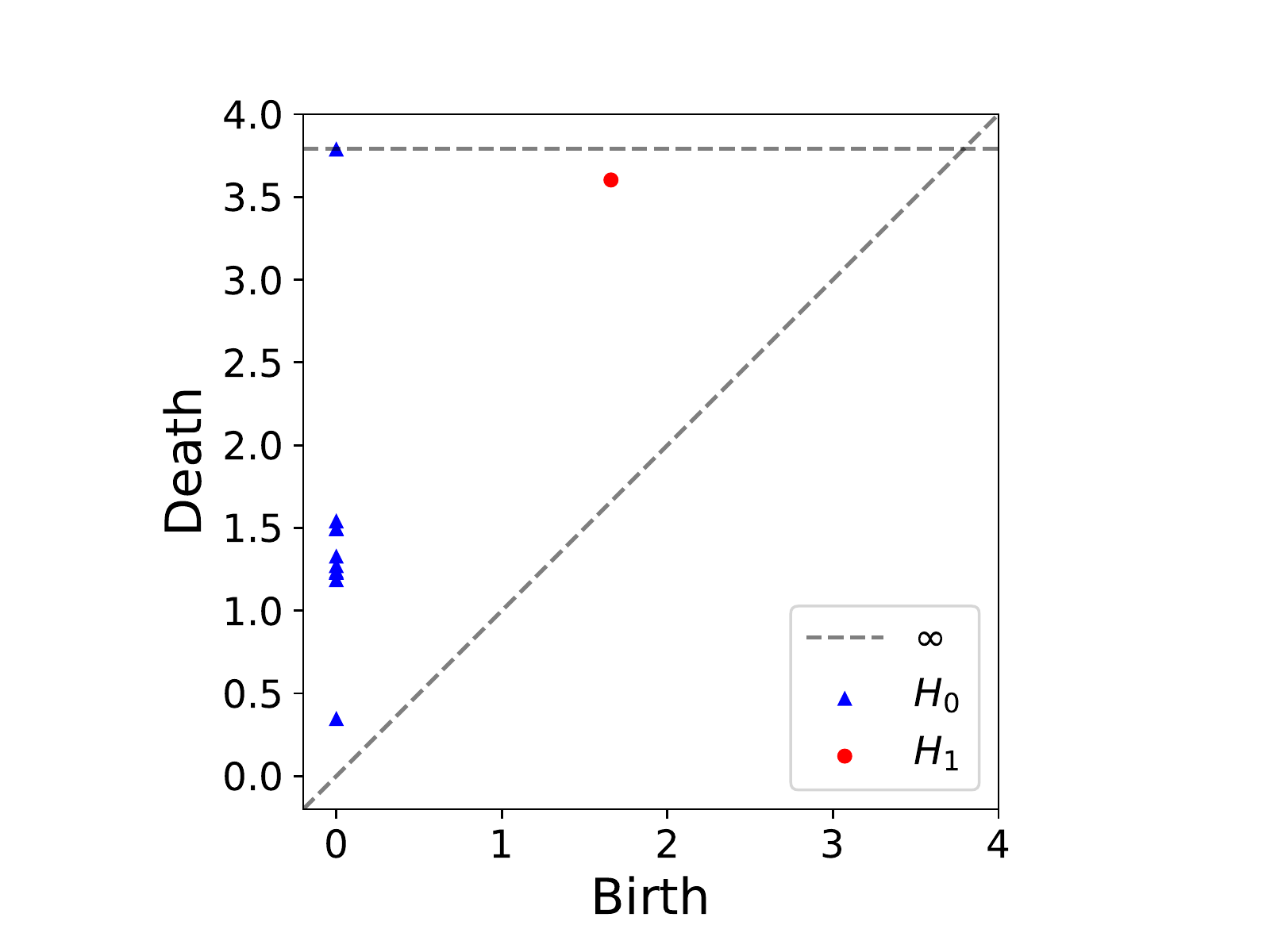}
	\caption{Persistent homology plot of the illustrative point cloud. The x-axis shows the birth and the y-axis the death of the homology groups}
	\label{fig:LVD_Method_Illustration_Persistent_Plot}
\end{figure}

In persistent homology, we are interested in so-called topological invariants, i.e., properties that are invariant under homeomorphisms. 
The topological invariants of interest are homology groups, i.e., H$_k$ of dimension $k$, with $\beta_k = \dim (H_k)$ being the Betti numbers~\citep{binchi2014jholes,chung2015persistent}. 
``Informally, $\beta_0$ is the number of connected components, $\beta_1$ is the number of two-dimensional holes or ``handles'' and $\beta_2$ is the number of three-dimensional holes or ``voids'' etc.''~\citep{binchi2014jholes}. \newline \indent
The topological invariants are computed by representing the original dataset, i.e., a point cloud, as a simplicial complex through a simplicial filtration. 
We utilize the common Vietoris-Rips filtration, where a n-simplex in the simplicial complex is formed if and only if the pairwise distance between all points in the n-simplex is at most $\epsilon$.  
At the bottom of Figure~\ref{fig:LVD_Method_Illustration_Rips_filtration}, we show a series of simplicial complexes for an illustrative point cloud. \newline \indent
Persistent homology studies topological invariants that persist over multiple length scales ($\epsilon$) in the data~\citep{chambers2018persistent,otter2017roadmap,xia2018persistent,xia2019persistent}.
In other words, we examine the lifespan of topological invariants by increasing $\epsilon$ incrementally and constructing simplicial complexes. 
At the bottom of Figure~\ref{fig:LVD_Method_Illustration_Rips_filtration}, we can observe that 
$\beta_0=10$ connected components (H$_0$) exist at $\epsilon_1$, 
$\beta_0=5$ connected components exist at $\epsilon_2$, 
$\beta_0=1$ connected component and $\beta_1=1$ two-dimensional hole (H$_1$) exist at $\epsilon_3$, and 
$\beta_0=1$ connected component exist at $\epsilon_4$. \newline \indent
The results of the persistent homology can be depicted in barcode diagrams or persistent diagrams.
We use the more common persistent diagrams in this work. 
The coordinates of birth and death of the homology groups in the example are shown in the persistent diagram in Figure~\ref{fig:LVD_Method_Illustration_Persistent_Plot}. 
The x-axis represents the $\epsilon_{\text{birth}}$ while the y-axis the $\epsilon_{\text{death}}$ distance of H$_0$ and H$_1$ homology groups.
Features with long lifespan correspond to points far from the diagonal~\citep{wasserman2018topological}.
The blue triangle at the bottom left corner of the plot corresponds to the merge of two very close data points at small $\epsilon$. 
The blue triangles with $\epsilon_{\text{death}}$ between 1 and 1.5 in Figure~\ref{fig:LVD_Method_Illustration_Persistent_Plot} resemble the merge of connected components between $\epsilon_2$ and $\epsilon_3$ in Figure~\ref{fig:LVD_Method_Illustration_Rips_filtration}, describing the decrease in Betti number $\beta_0$ from 5 to 1. 
The highest blue triangle illustrates that one connected component exists until infinite. 
The red circle represents homology group H$_1$ and demonstrates the birth and death of the two-dimensional hole which is formed around $\epsilon_3$ and dies before $\epsilon_4$ in Figure~\ref{fig:LVD_Method_Illustration_Rips_filtration}. \newline \indent
In this example, the persistent diagram shows that a hole exist providing useful insight to guide the decision process for model selection.
For example, the $\epsilon_{\text{death}}$ of the H$_1$ components provide information about the data density.
In the example, the maximal distance in the last H$_1$ component is at most 1.5.
This is significantly smaller than the $\epsilon_{\text{death}}$ of the H$_1$ hole.
In other words, the life span of the hole is characteristic for the dataset. 
\subsection{Learn validity domain using one-class support vector machines}
\label{sec:Lear_validity_domain_Background_One_Class_SVM}
We model the validity domain using the convex hull or one-class SVM approach.
The one-class SVMs are trained using the open-source python implementation in Scikit-learn~\citep{pedregosa2011scikit} and the convex hulls are identified using the open-source implementation in SciPy~\citep{SciPy}.
The details of the one-class SVM are described in the following. \newline \indent
SVMs are a popular method for binary classification~\citep{cortes1995support} and regression~\citep{smola2004tutorial}.
One-class SVMs are a modification of these classical SVMs~\citep{scholkopf2000support} (c.f. \citet{Tax2001} on similarity to SVDD). 
The goal is to learn a boundary of a given set of training points ${X} = \{\hat{\boldsymbol{x}}^{(1)}, ...,\hat{\boldsymbol{x}}^{(i)},..., \hat{\boldsymbol{x}}^{(N)}\}$ with $\hat{\boldsymbol{x}}^{(i)} \in \mathbb{R}^{D}$.
Similar to classical SVMs, the data is mapped to  high-dimensional feature space by $\phi : \mathbb{R}^{D} \mapsto \mathbb{R}^{d}$ with $d >> D$ and later solved in the dual formulation using the \textit{kernel trick}~\citep{scholkopf2001kernel}.
In the feature space, a maximum-margin hyperplane is found that separates the data from the origin by solving:
\begin{align}
&\min_{\boldsymbol{w} \in \mathbb{R}^{d}, \xi_i \in \mathbb{R}, \rho \in \mathbb{R}} 	&&\frac{1}{2}\boldsymbol{w}^T\boldsymbol{w} -\rho + \frac{1}{\nu N}\sum_{i=1}^{N} \xi_i	\label{eqn:one_class_SVM_primal_1}, \\
&\text{s.t} &&{\boldsymbol{w}^T\phi(\hat{\boldsymbol{x}}^{(i)})\geq \rho- \xi_i} & \forall i =\{1,..,i,..N\} \label{eqn:one_class_SVM_primal_2},\\
&			&& \xi_i \ge 0 	& \forall i =\{1,..,i,..N\} \label{eqn:one_class_SVM_primal_3},
\end{align}
where $\nu \in (0,1) $ is a regularization hyperparameter, $\xi_i \in \mathbb{R}$  are slack variables, and $\boldsymbol{w}$ and $\rho$ are the parameters of the hyperplane in high-dimensional feature space.  
\cite{scholkopf2000support} show that $\nu$ is an upper bound on the fraction of outliers and a lower bound on the fraction of support vectors in the training  set, which is known as the  $\nu$-property.
The decision function $f_{\text{DF-P}}(\boldsymbol{x}) = \boldsymbol{w}^T\phi(\boldsymbol{x})-\rho$ is positive if a candidate point $\boldsymbol{x}$ is classified to be within the training data domain and negative if not.
The dual of \eqref{eqn:one_class_SVM_primal_1}-\eqref{eqn:one_class_SVM_primal_3} is the quadratic program:
\begin{align}
&\min_{\alpha_i \in [0,\frac{1}{\nu N}]} 
&&\frac{1}{2} \sum_{i=1}^{N} \sum_{j=1}^{N} \alpha_i  \alpha_j K(\hat{\boldsymbol{x}}^{(i)},\hat{\boldsymbol{x}}^{(j)}), \\
&\text{s.t} &&{\sum_{i=1}^{N} \alpha_i = 1},
\end{align}
where $\alpha_i$ are dual variables and $K(\cdot,\cdot)$ is a kernel function. 
Often, the radial basis kernel function  $K(\boldsymbol{x},\boldsymbol{y}) =  \exp(-\gamma\|\boldsymbol{x}-\boldsymbol{y}\|^2)$ with hyperparameter $\gamma$ is used since it has been shown that it is best able to model the most complex boundaries~\citep{Tax2001}.
It holds that $\alpha_i =  0 $ for training samples inside the learned boundary and $\alpha_i > 0 $  for samples on or outside the boundaries. 
Samples for which $\alpha_i > 0$ are called support vectors.
The decision function in the dual variables is given by $f_{\text{DF-D}}(\boldsymbol{x}) = \sum_{i \in I_{\text{sv}}}  \alpha_i K(\hat{\boldsymbol{x}}^{(i)}, \boldsymbol{x}) -\rho$, where $I_{\text{sv}}$ denotes the indexes of the support vectors in the training data (i.e., the data points with corresponding $\alpha_i > 0$).
To obey validity limits of data-driven models in an optimization problem, the following inequality has to hold:
\begin{equation}
f_{\text{DF-D}}(\boldsymbol{x})  \geq 0 \label{eq:decision_function}.
\end{equation}
The parameter $\nu$ can be estimated from an outlier fraction by using the aforementioned $\nu$-property. 
This makes this method more tolerant to outliers in the training data~\citep{pimentel2014review}. 
The hyperparameter $\gamma$ controls the model complexity when using the radial basis kernel.
If a large $\gamma$ is used, all samples are mapped to a small region in the feature space and the one-class SVM cannot distinguish between the samples well. 
In other words, the model lacks complexity. 
If $\gamma$ is small, pairs of samples become orthogonal in the feature space. 
This leads to overfitting and a high number of support vectors. 
A common approach to identify an appropriate $\gamma$ is to gradually decrease $\gamma$ until the number of support vectors does not decrease much~(e.g., \cite{dreiseitl2010outlier}). 
However, automatically selecting an appropriate $\gamma$ is challenging (e.g.,\citep{evangelista2007some,xiao2014parameter,xiao2014two}).
\subsection{Optimization with classifier as constraint}
\label{sec:Lear_validity_domain_Method_Optimization}
We consider a global optimization problem where a classifier is used to obey validity limits of a data-driven model. 
In most cases, the inputs of the classifier model correspond to the degrees of freedom $\boldsymbol{x}$ of the optimization problem. 
The classifier can determine if a given $\boldsymbol{x}$ is feasible or infeasible by evaluating its decision function $f_{\text{DF-D}}(\cdot)$.
To obey validity limits, Inequality~\eqref{eq:decision_function} has to hold. 
Although the decision function is an explicit function, there exist different ways to formulate it in optimization problems. 
These problem formulations are equivalent as they have the same solution, but they can have a large impact on the computational performance of global optimization. \newline \indent
In the FS formulation, a set of nonlinear equations is provided as equality constraints while the dependent (or intermediate) variables are optimization variables.
Note that there exist multiple valid FS formulations depending on the equality constraints and optimization variables provided to the solver. 
One representative FS formulation for optimization with one-class SVMs embedded is shown in the following: 
\begin{align}
&\min_{\boldsymbol{x} \in \mathbb{R}^D,\newline
	z_{\text{obj}} \in \mathbb{R},
	d_i  \in \mathbb{R}, \newline
	k_i \in \mathbb{R},
	\boldsymbol{z}_{\text{dd}} \in  \mathbb{R}^{Z}
} &&z_{\text{obj}} \label{eq:LFR_FS_obj},\\
&\text{s.t} && z_{\text{obj}} = f_{\text{obj}}(\boldsymbol{x},\boldsymbol{z}_{\text{dd}}) \label{eq:LFR_FS_con1},\\
&&& \boldsymbol{h}_{\text{dd}}(\boldsymbol{x},\boldsymbol{z}_{\text{dd}})=\boldsymbol{0} \label{eq:LFR_FS_con2},\\
&&&\sum_{i \in I_{\text{sv}}}\alpha_i k_i \geq \rho \label{eq:LFR_FS_con3},\\
&&&{k_i}{= \exp{(-\gamma \cdot d_i)} } &{\forall i \in I_{\text{sv}}}\label{eq:LFR_FS_con4},\\
&&&{d_i}{ = \|\hat{\boldsymbol{x}}^{(i)} - \boldsymbol{x}\|^2 } &{\forall i \in I_{\text{sv}}} \label{eq:LFR_FS_con5}.
\end{align}
Herein, Equation~\eqref{eq:LFR_FS_obj} minimizes the objective function value $z_{\text{obj}}$ that is given by Equation~\eqref{eq:LFR_FS_con1}.
Note that the objective depends on the variables of the data driven model $\boldsymbol{z}_{\text{dd}}$ that are given by the solution of Equations~\eqref{eq:LFR_FS_con2}.
The decision function of the one-class SVM is given by the inequality constraint~\eqref{eq:LFR_FS_con3} while its intermediate variables are given by the solution of Equations~\eqref{eq:LFR_FS_con4},\eqref{eq:LFR_FS_con5}. 
This FS formulation has in total $D+2\cdot \vert I_{\text{sv}}\vert +\dim(\boldsymbol{z}_{\text{dd}})+ 1$ optimization variables, $2\cdot \vert I_{\text{sv}}\vert +\dim(\boldsymbol{z}_{\text{dd}})+ 1$ equality constraints, and one inequality constraint. \newline \indent
The equality constraints of the one-class SVM can be solved explicitly for the intermediate variables. 
Thus, we can directly formulate a RS formulation of the optimization problem~(c.f.~\citep{Bongartz.2017}):
\begin{align}
&\min_{\boldsymbol{x} \in \mathbb{R}^D}&&f_{\text{RS}}(\boldsymbol{x}) \label{eq:LFR_RS_obj}, \\
&\text{s.t} &&f_{\text{DF-D}}(\boldsymbol{x})  \geq 0  \label{eq:LFR_RS_con1}.
\end{align}
Herein, $f_{\text{RS}}(\cdot)$ is the RS formulation of the data-driven model and objective function.
Thus, Equation~\eqref{eq:LFR_RS_obj} results from sequential substitutions of Equations~\eqref{eq:LFR_FS_obj}-\eqref{eq:LFR_FS_con2}.
This is possible as most data-driven models such as ANNs or GPs are explicit functions (c.f. \citet{Schweidtmann2019detGlobalANN,schweidtmann2020globalGP}) and as the objective function is a function of the the degrees of freedom and the predictions of the data-driven model.
Similarly, Equation~\eqref{eq:LFR_RS_con1} results from the substitution of Equations~\eqref{eq:LFR_FS_con3}-\eqref{eq:LFR_FS_con5}.
The RS formulation has only $D$ optimization variables and one inequality constraint. \newline \indent
The convex hull of a point cloud with a finite number of points can be formulated as a set of linear inequality constraints ${X}_\text{convHull} = \{\boldsymbol{x} \in \mathbb{R}^D~|~\boldsymbol{A}\boldsymbol{x}+\boldsymbol{b} \leq \boldsymbol{0}\} $. 
Assuming that the convex hull has $f$ facets, the matrix $\boldsymbol{A} \in \mathbb{R}^{f\times D}$ and the vector $\boldsymbol{b} \in \mathbb{R}^f$~\citep{kahrs2007validity}. 
Thus, the FS and RS formulation of the convex hull are identical. 
Note that the data-driven model can still be formulated in the RS and FS formulation when using the linear convex hull constraints. \newline \indent
The RS formulation has three major advantages for global optimization:
First, the problem formulation has a direct influence on the variables to be branched on. 
In the RS, the B\&B solver branches only on the degrees of freedom $\boldsymbol{x}$.
In the FS, the B\&B solver branches on the degrees of freedom and also on the intermediate variables.
This is undesirable given the exponential worst-case runtime of global optimization methods.
Note that this issue can also be mitigated by selective branching \citep{Epperly.1997}.
Second, the size of the subproblems that are solved during optimization is affected by the problem formulation and the method for constructing relaxations. 
Our previous work shows that a combination of McCormick relaxations and RS formulation can reduce the time to solve an iteration of the B\&B solver significantly~\citep{schweidtmann2020globalGP,bongartz2020_diss}. 
Third, global optimization solvers usually require bounds on all optimization variables. 
Often, meaningful bounds are known for degrees of freedom but bounds on intermediate variables can be difficult to determine. 
Note that this problem is mitigated by some state-of-the-art solvers through automatic bound tightening techniques.  \newline \indent
The vast majority of previous literature approaches formulate global optimization in the FS because they frequently use modeling environments such as GAMS that essentially require an equation-oriented modeling approach. 
Recently, \citet{hart2017pyomo} developed the Python-based optimization tool Pyomo which allows modeling, implementation of own solvers, and provides access to multiple solvers. 
Pyomo allows both FS and RS and recently \citet{hullen2019managing} demonstrated RS optimization of ANNs in BARON through Pyomo. 
However, BARON relies on the auxiliary variable method for relaxations which results in larger subproblems~\citep{schweidtmann2020globalGP}. 
Thus, a RS formulation in BARON does not take full advantage of the RS formulation. 
In contrast, our open-source solver MAiNGO~\citep{MAiNGO} relies on McCormick relaxations in the space of the original variables utilizing the library MC++~\citep{Mitsos.2009,Chachuat.2015}. 
Another open-source solver that allows for McCormick relaxations is called EAGO has been released by \citet{wilhelm2020eago}.  \newline \indent
The optimization problems in this work are implemented in MeLOn~\citep{MeLOn_Git} and solved by MAiNGO~\citep{MAiNGO}. 
For comparison, the problems are additionally exported to GAMS and solved by the commercial solver BARON~\citep{Tawarmalani.2005}. 
We provide the implementation of the one-class SVM in the open-source modeling toolbox MeLOn~\citep{MeLOn_Git}. \newline \indent
Tight convex and concave relaxations are highly desirable in global optimization. 
Therefore, we use the tightest possible relaxations, i.e., the envelopes, of the radial basis function kernel in our solver MAiNGO. 
Note that we derived these envelopes in our previous work~\citep{schweidtmann2020globalGP} as the squared exponential covariance function in Gaussian processes is equivalent to the radial basis function kernel in SVMs.
\section{Illustrative case studies}
As high dimensional problems are difficult to visualize, we consider eight two-dimensional data sets for illustration of the proposed method in a first step. 
Afterwards, we consider an engineering case study in Section~\ref{sec:Lear_validity_domain_engineering_case_study}. 
As shown in Figure~\ref{fig:LVD_Bounds_of_illustrative_example}, the illustrative examples cover a variety of relevant scenarios.
All data points are randomly generated within pre-specified bounds and perturbed by noise.
Thus, the data sets do not exhibit sharp boundaries, rather they also include noisy outlier data points. \newline \indent
In the next step, we evaluate an adapted peaks function, $f_{\text{Peaks}} : \mathbb{R}^2 \mapsto \mathbb{R}$, on all data sets with $f_{\text{Peaks}}(x_1,x_2) = 3(1-x_1)^2 \cdot \exp{(-x_1^2-(x_2+1)^2)} -10(\frac{x_1}{5}-x_1^3-x_2^5)\cdot \exp{(-x_1^2 -x_2^2)} - \frac{1}{3}\cdot \exp{(-(x_1+1)^2-x_2^2)}-1.3x_2$.
Then, we train individual ANNs on the eight data sets using Keras~\citep{chollet2015keras}.
All ANNs exhibit one input layer with 2 neurons, two hidden layers with six and eight neurons, respectively, and an output layer with one neuron.
The hidden layers use $\tanh$ activation and the output layers use linear activation.
For training, the inputs are scaled onto $[-1,1]$ and the outputs are scaled to zero mean and unit variance.
We further use a batch size of $128$ and an epoch limit of 4,000.
Note that we omit a hyperparameter study for the ANNs because ANN training is not the focus of this work. \newline \indent
All optimization problems are solved one core of an Intel Xeon CPU E5-2630 v3 (2.40GHz) with 192 GB RAM and Windows Server 2016 operating system.
We use a 0.001 relative and absolute optimality tolerance, a CPU time limit of 1,000 seconds, and default settings in BARON and MAiNGO.
\subsection{Topological data analysis}
\begin{figure*} 
	\begin{subfigure}[b]{0.50\textwidth}
		\centering
		\includegraphics[width = \textwidth]{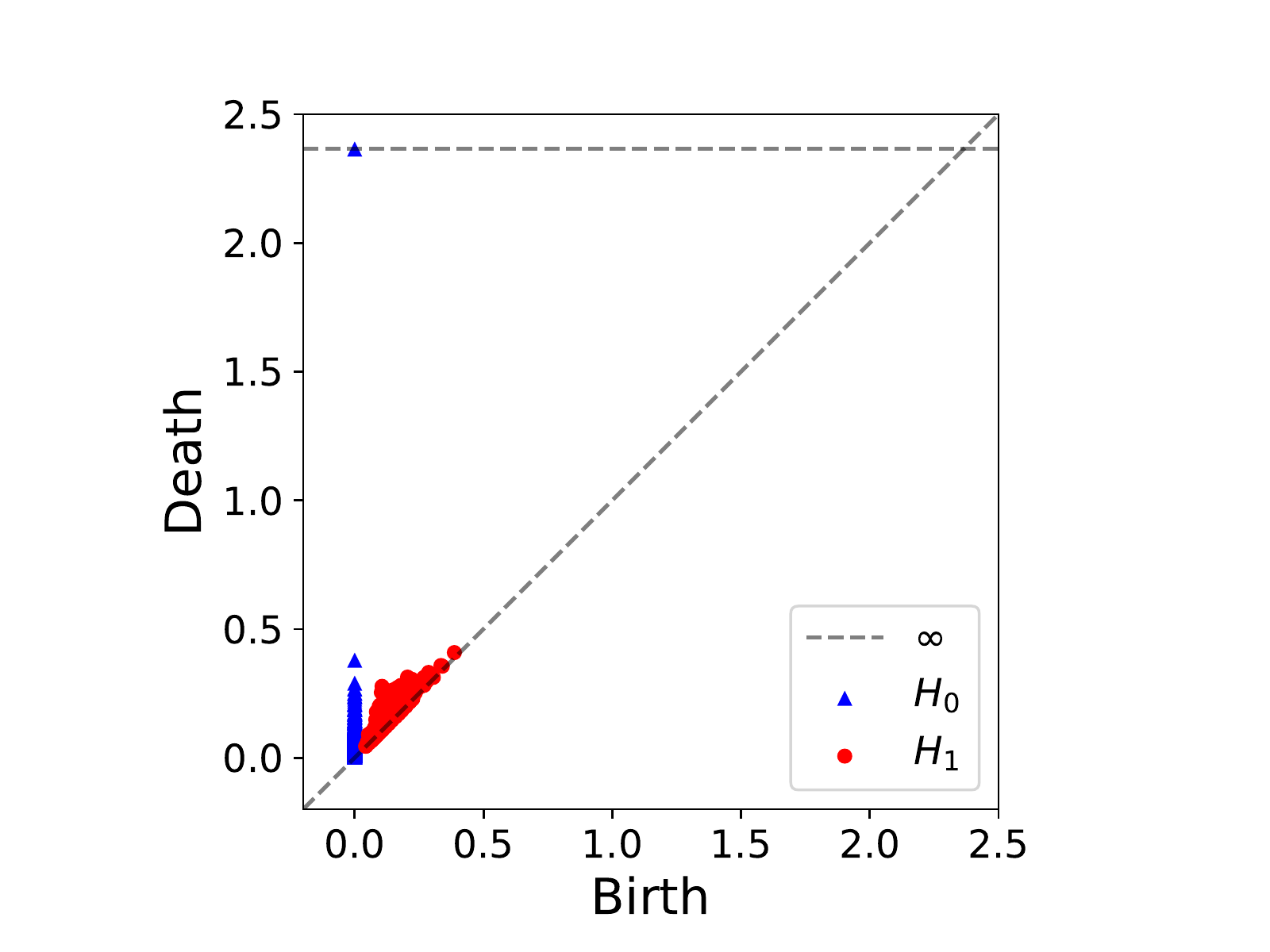}
		\caption{``Box'' case study}
		\label{fig:LVD_Persistent_Homology_Plots_box}
	\end{subfigure} \hfill
	\begin{subfigure}[b]{0.50\textwidth}
	\centering
	\includegraphics[width = \textwidth]{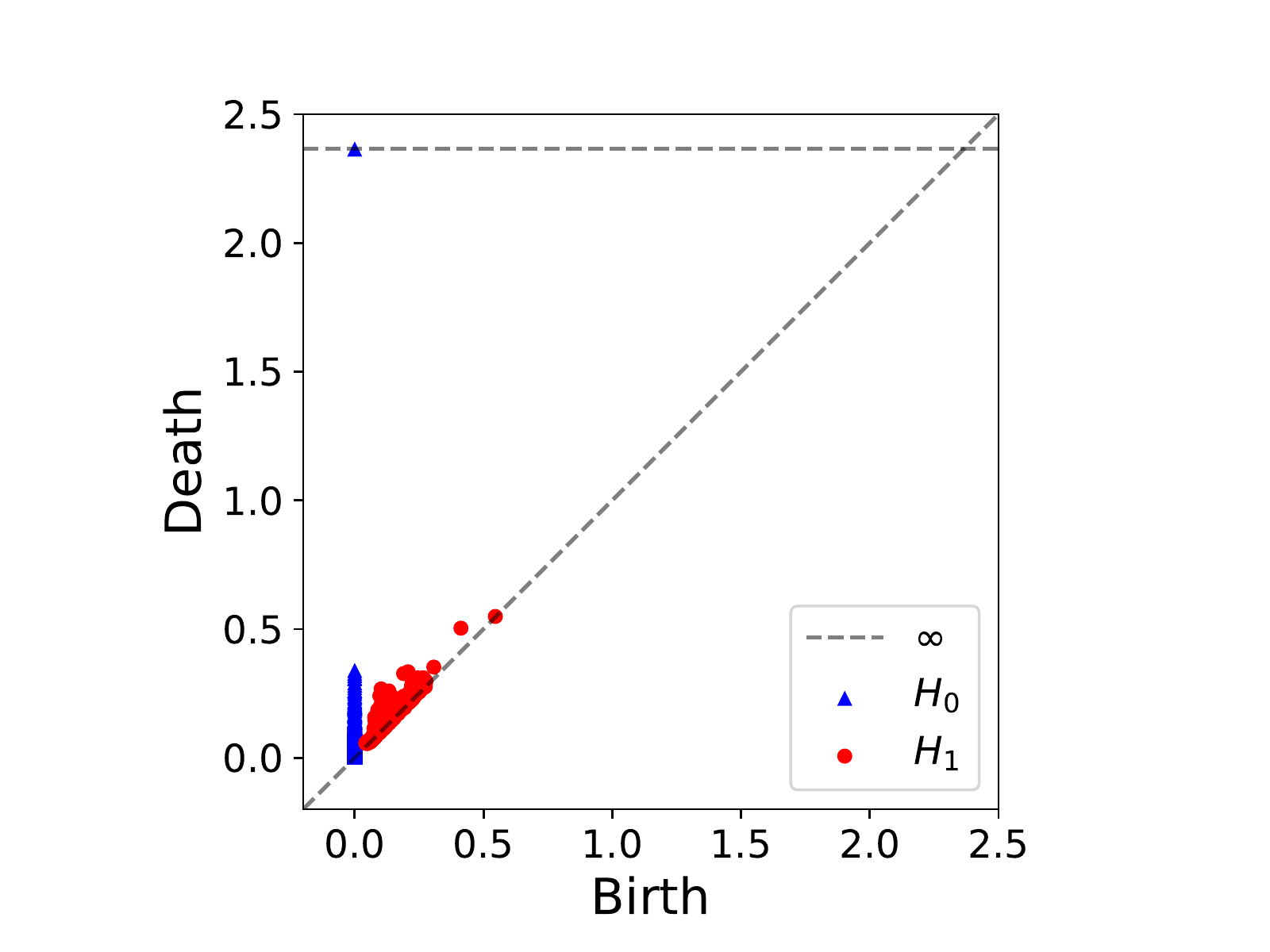}
	\caption{``Oval'' case study}
	\label{fig:LVD_Persistent_Homology_Plots_oval}
\end{subfigure} \\
	\begin{subfigure}[b]{0.50\textwidth}
	\centering
	\includegraphics[width = \textwidth]{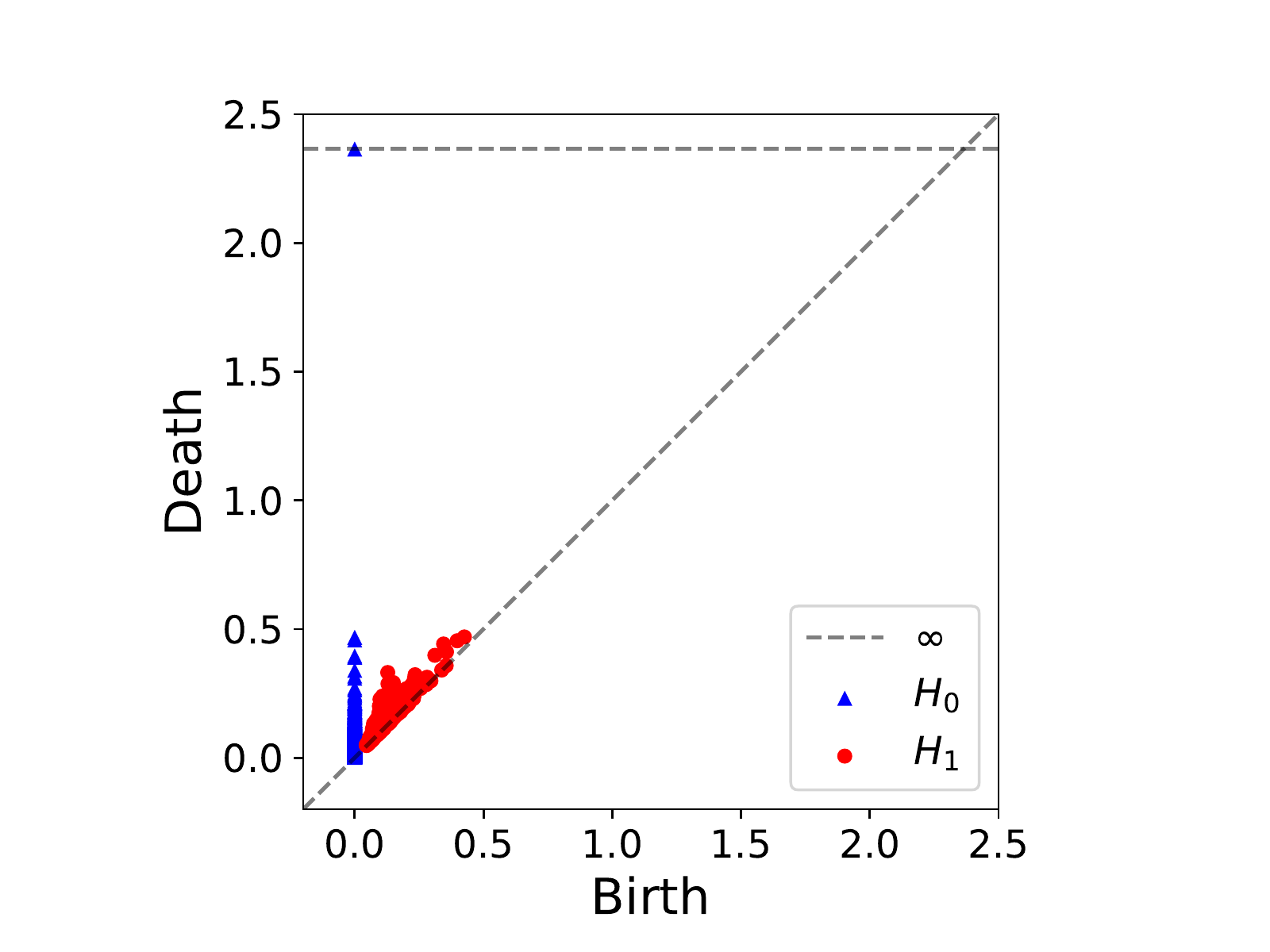}
	\caption{``Box2'' case study}
	\label{fig:LVD_Persistent_Homology_Plots_box2}
\end{subfigure} \hfill
	\begin{subfigure}[b]{0.50\textwidth}
	\centering
	\includegraphics[width = \textwidth]{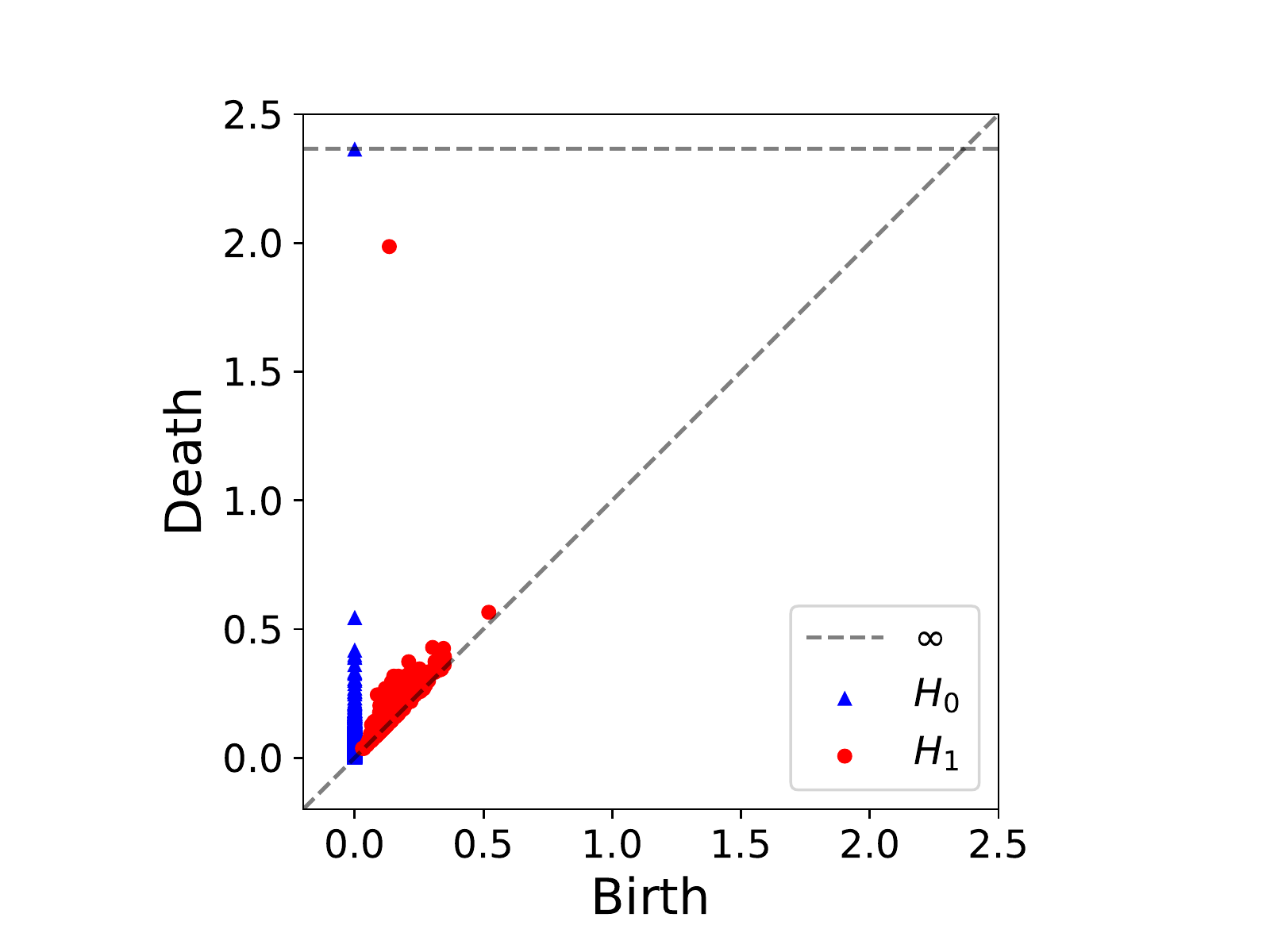}
	\caption{``Box w/ hole'' case study}
	\label{fig:LVD_Persistent_Homology_Plots_box_with_hole}
\end{subfigure} \\
	\begin{subfigure}[b]{0.50\textwidth}
	\centering
	\includegraphics[width = \textwidth]{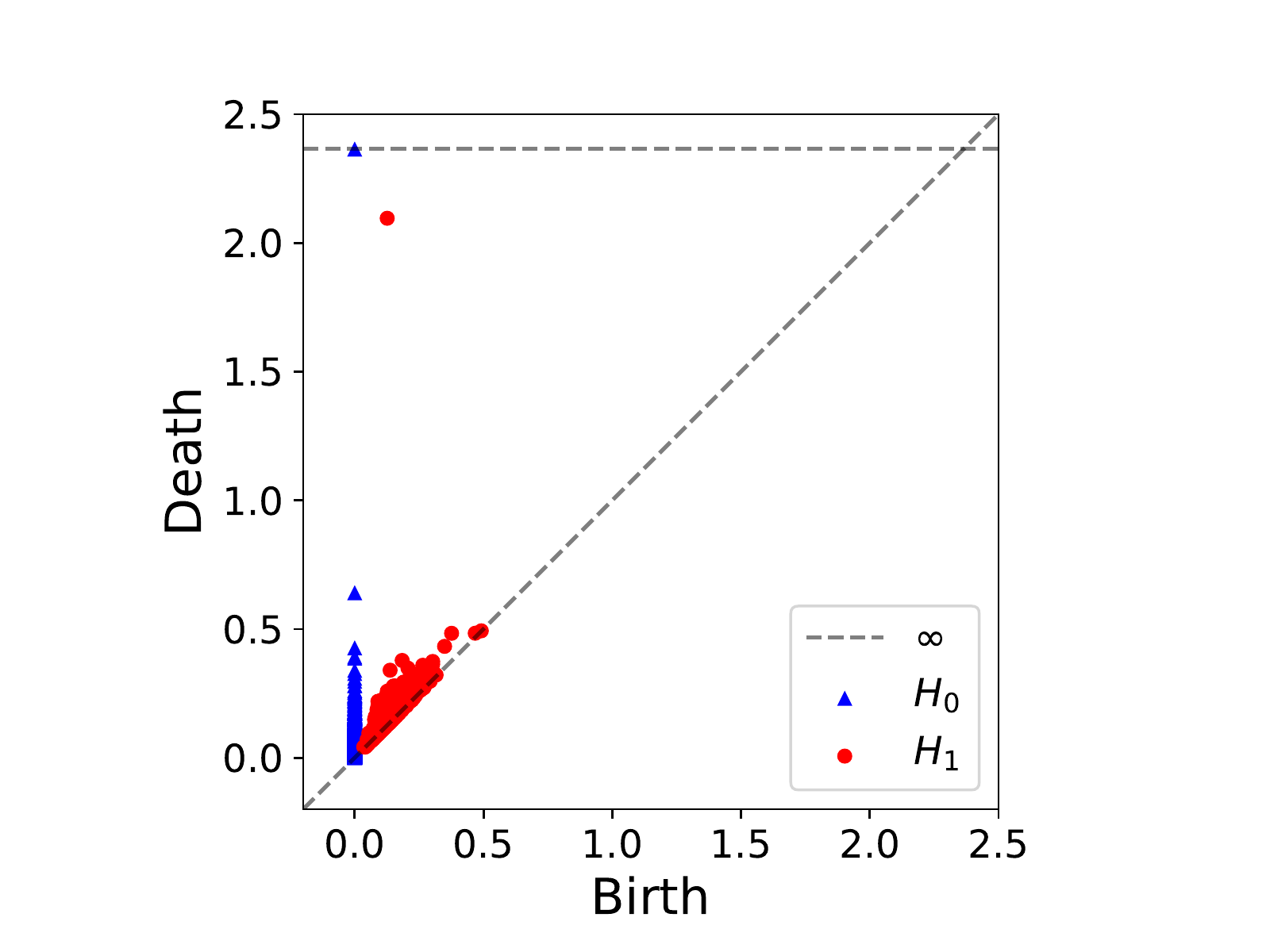}
	\caption{``Circle w/ hole'' case study}
	\label{fig:LVD_Persistent_Homology_Plots_circle_with_hole}
\end{subfigure} \hfill
	\begin{subfigure}[b]{0.50\textwidth}
	\centering
	\includegraphics[width = \textwidth]{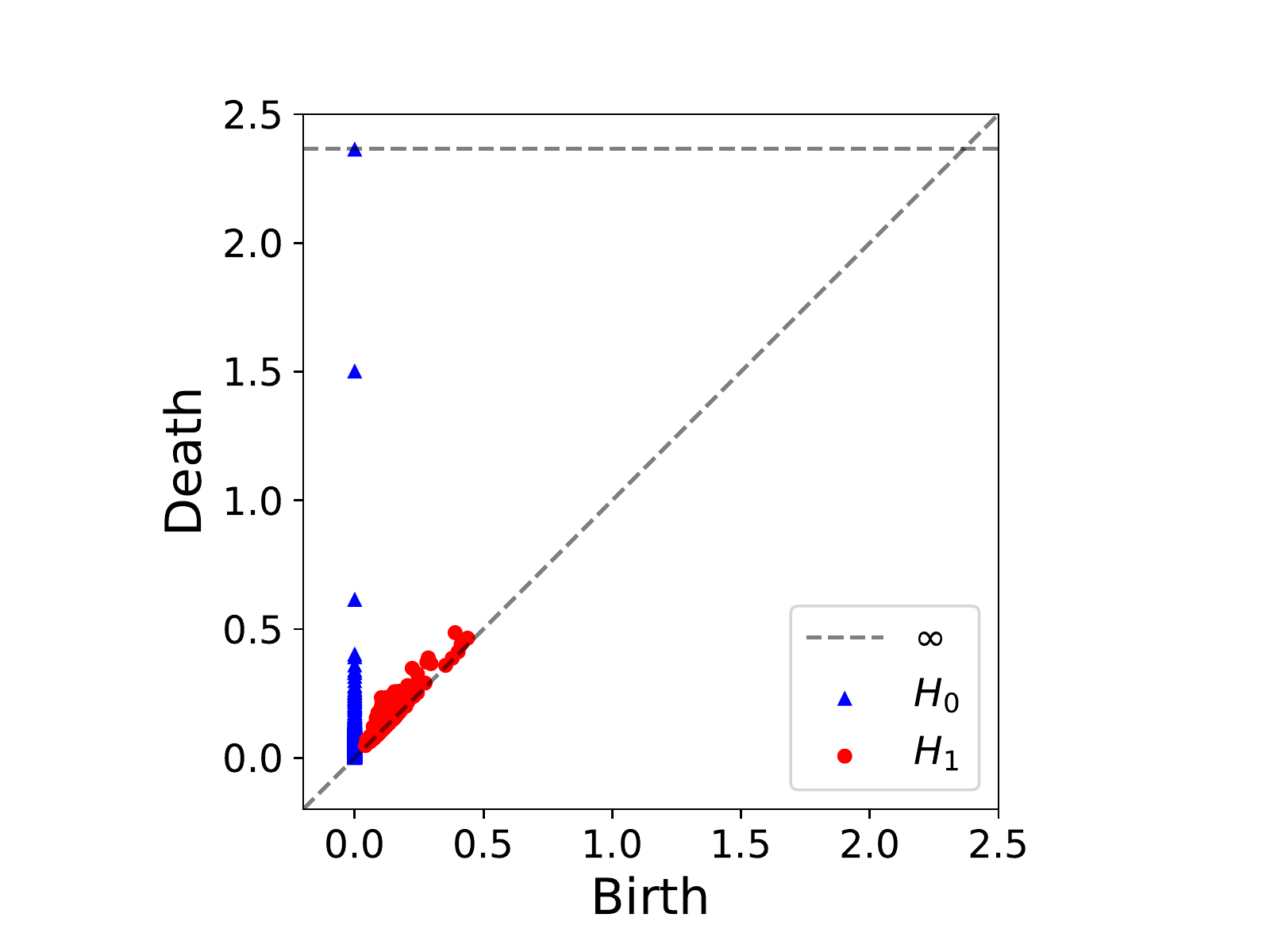}
	\caption{``Two circles'' case study}
	\label{fig:LVD_Persistent_Homology_Plots_two_circles}
\end{subfigure} \\
	\begin{subfigure}[b]{0.50\textwidth}
	\centering
	\includegraphics[width = \textwidth]{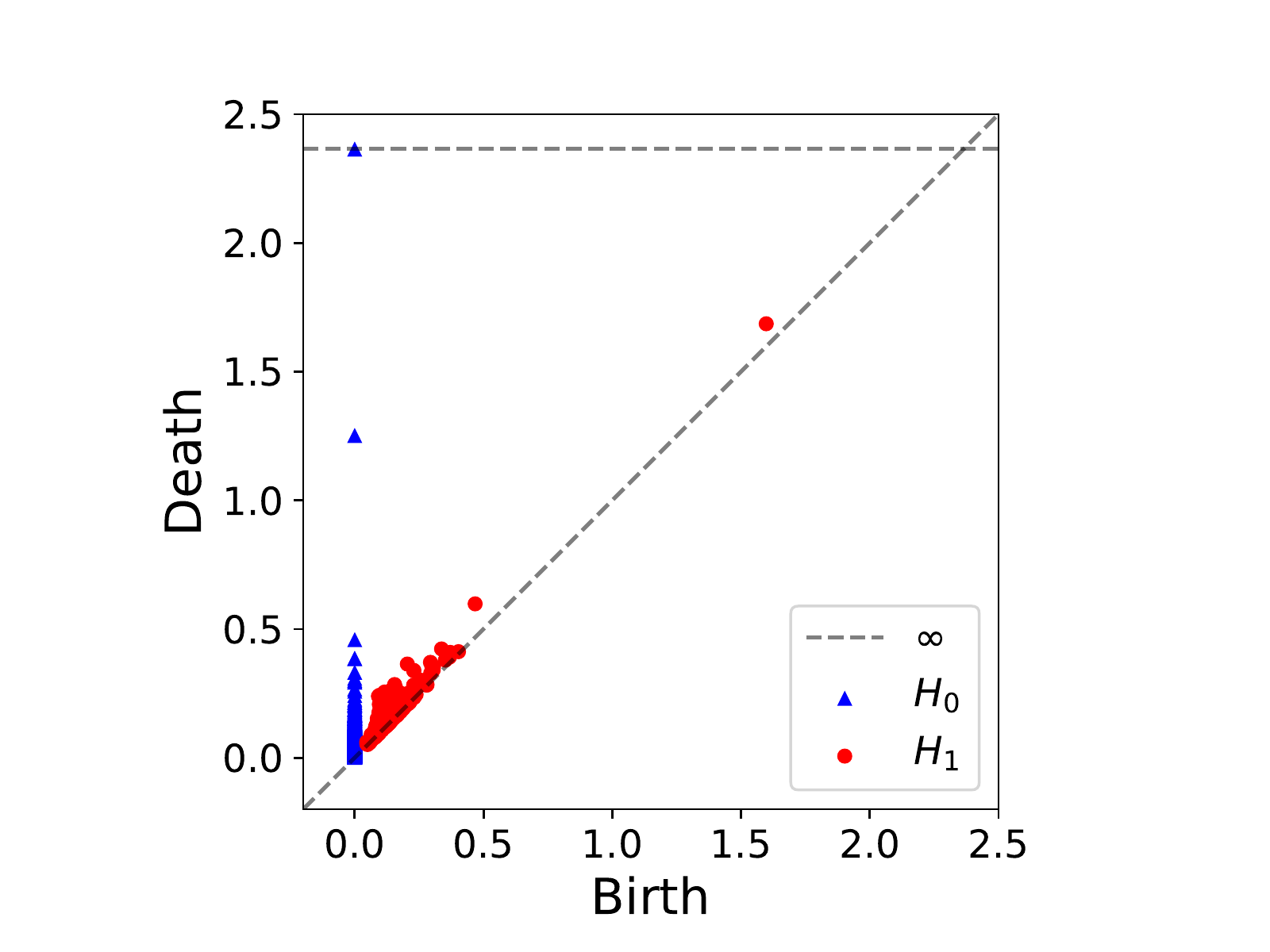}
	\caption{``Two ovals'' case study}
	\label{fig:LVD_Persistent_Homology_Plots_two_ovals}
\end{subfigure} \hfill
	\begin{subfigure}[b]{0.50\textwidth}
	\centering
	\includegraphics[width = \textwidth]{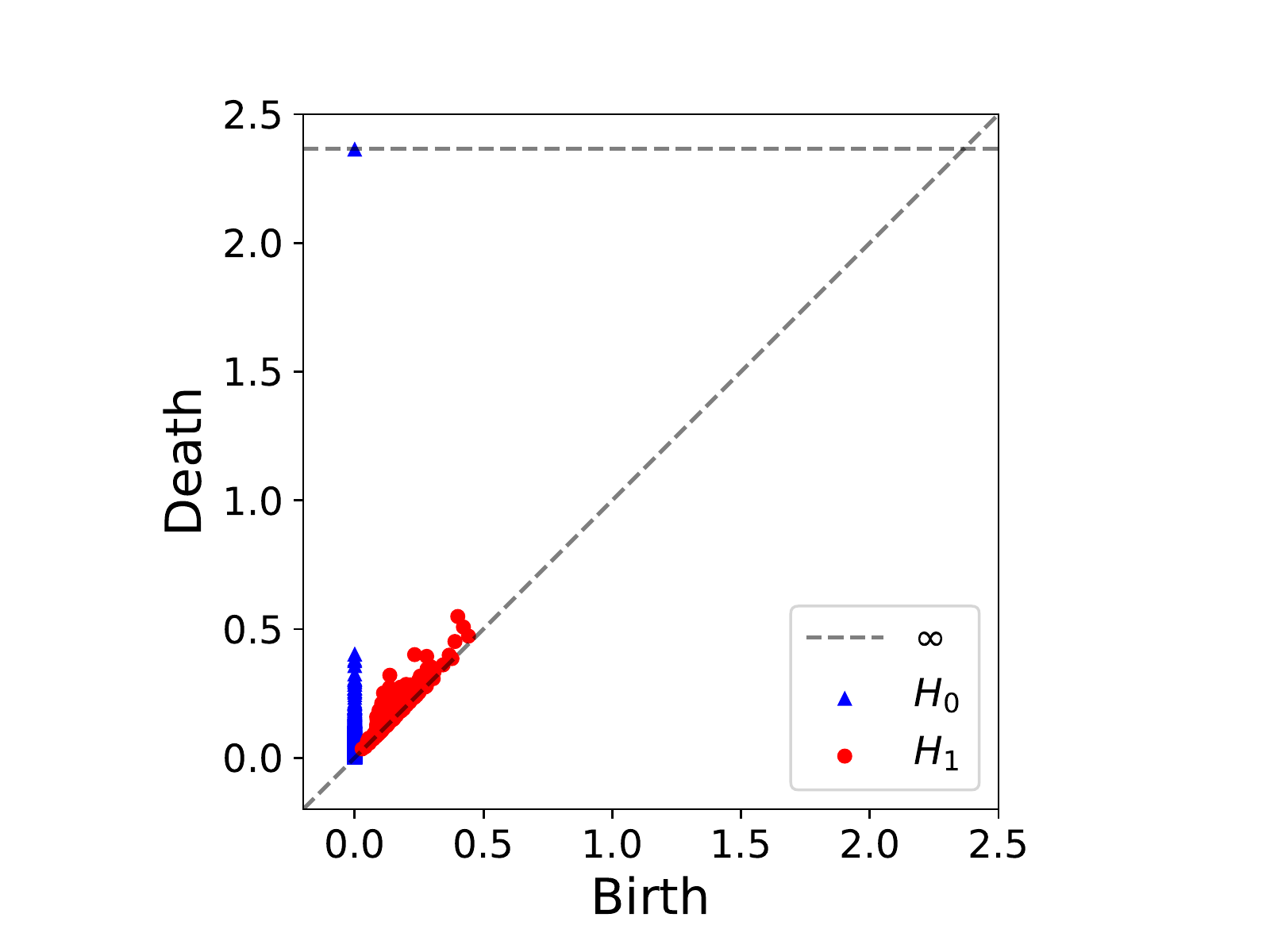}
	\caption{``Banana'' case study}
	\label{fig:LVD_Persistent_Homology_Plots_banana}
\end{subfigure}
	\caption{Comparison of the persistent homology plots of the eight illustrative data sets}
	\label{fig:LVD_Persistent_Homology_Plots}
\end{figure*}

The persistent diagrams of the eight input data sets are shown in Figure~\ref{fig:LVD_Persistent_Homology_Plots}. 
The Subfigures~\ref{fig:LVD_Persistent_Homology_Plots_box_with_hole} and \ref{fig:LVD_Persistent_Homology_Plots_circle_with_hole} show a H$_1$ component with a large life span each.
These correspond to the holes in the respective data sets ``box w/ hole'' and ``circle w/ hole''.
Moreover, the corresponding $\epsilon_{\text{death}}$ provide information about the diameter of the holes.
Recall that H$_1$ components that are close to the diagonal have a short life span and are therefore not relevant for this analysis. \newline \indent
The persistent diagrams also show the existence of disjunct clusters in the data sets ``two circles'' and ``two ovals''.
In Subfigures~\ref{fig:LVD_Persistent_Homology_Plots_two_circles} and \ref{fig:LVD_Persistent_Homology_Plots_two_ovals}, the H$_0$ components with a high $\epsilon_{\text{death}}$ indicate disjunct clusters and are a measure for the distance between them.
Note that the $H_0$ components at the infinity line persist when $\epsilon$ goes to infinity and do not die.
They correspond to the $H_0$ components that include all data points. \newline \indent
The persistent diagrams show no distinct differences between the ``box'', ``oval'', ``box2'', and ``banana'' case studies.
This illustrates that the method cannot distinguish between convex and nonconvex data sets in general.
Therefore, the persistent diagrams cannot ensure that the convex hull is sufficient to describe the validity domain.
Rather, it can only identify some cases where the convex hull is not sufficient.
\subsection{Validity domain modeling}
\begin{figure}
\begin{subfigure}[b]{\textwidth}
	\centering
	\includegraphics[width = \textwidth]{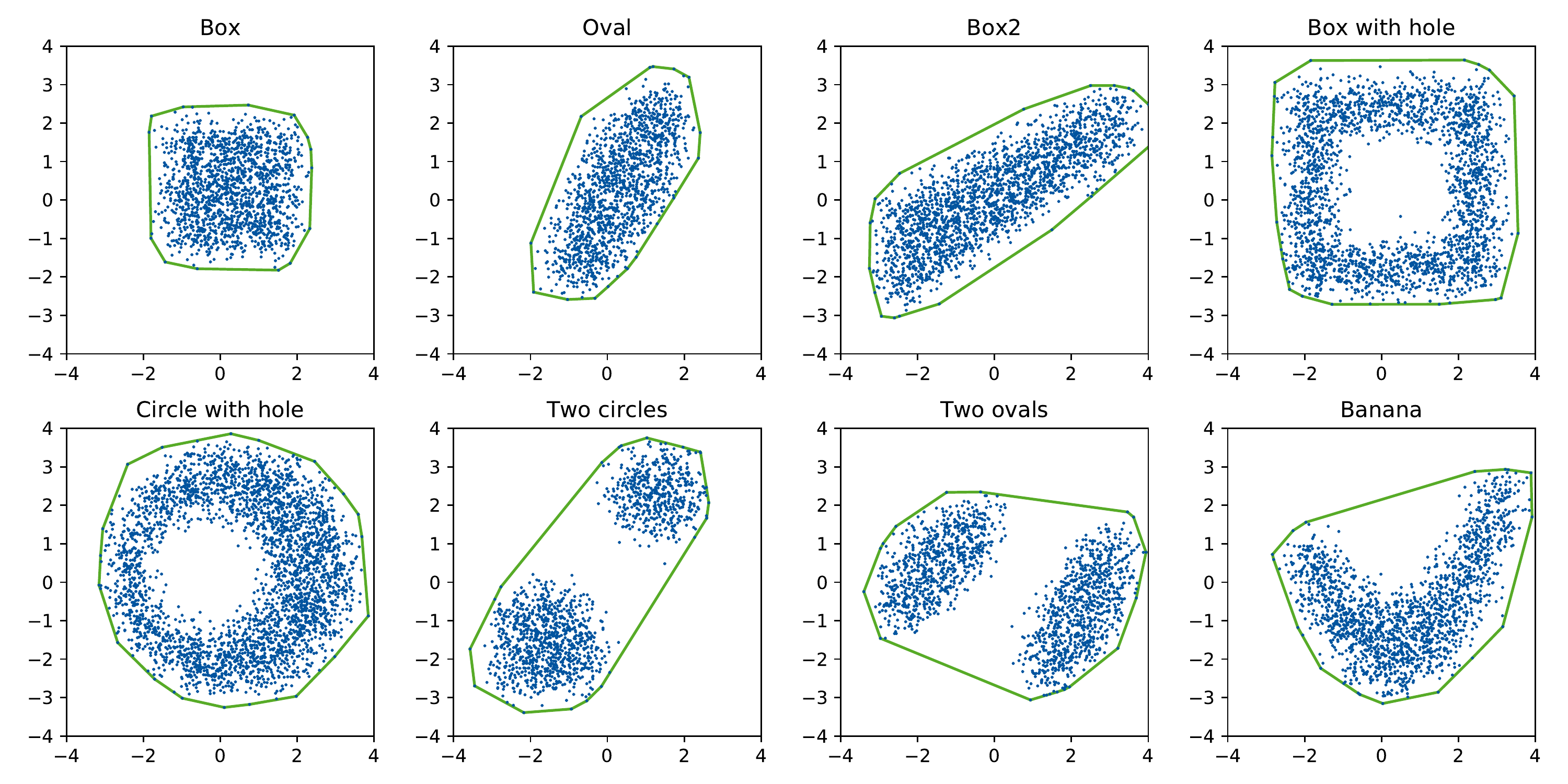}
	\caption{Convex hulls of the eight illustrative data sets}
	\label{fig:LVD_Convex_hull_of_illustrative_example}
\end{subfigure} \\
\begin{subfigure}[b]{\textwidth}
	\centering
	\includegraphics[width = \textwidth]{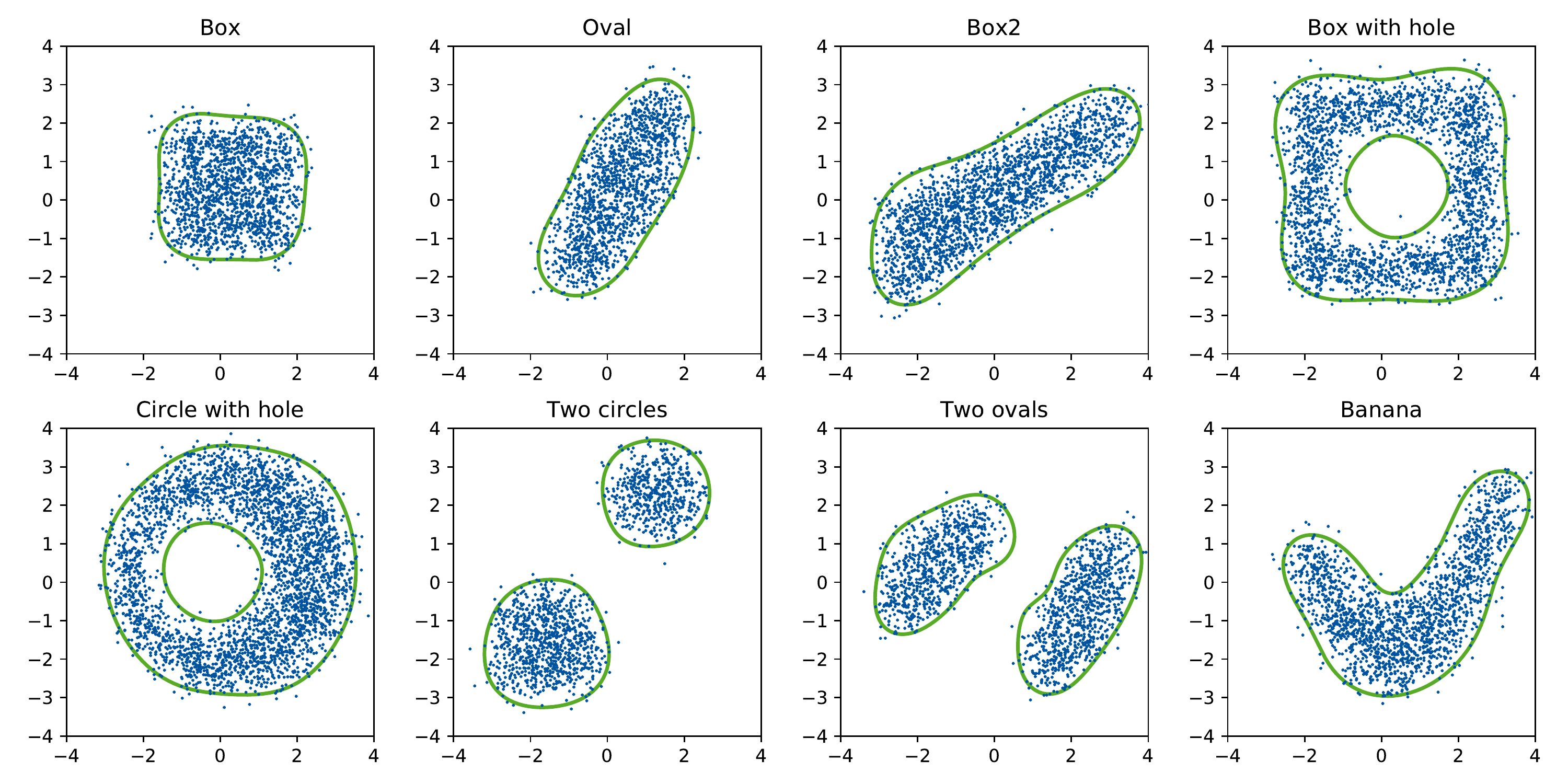}
	\caption{One-class SVM boundaries of the eight illustrative data sets}
	\label{fig:LVD_SVM_bounds_of_illustrative_example}
\end{subfigure}
	\caption{Comparison of the convex hull and the boundaries learned by the one-class SVMs}
\label{fig:LVD_Bounds_of_illustrative_example}
\end{figure}

In order to compare the one-class SVM and the convex hull, we model the input data of the eight case studies with both methods.
We employ the common radial basis function kernel for the one-class SVM. 
We set the hyperparameter $\nu$ to a low value $0.03$ because there are only a few outliers through noise in the data. 
The hyperparameter $\gamma$ is identified using the incremental approach described in Section~\ref{sec:Lear_validity_domain_Background_One_Class_SVM}.  
The selected $\gamma$ values are summarized in Table~\ref{tab:Lear_validity_domain_GAMMA}. \newline \indent
\begin{table}
	\centering
	\begin{tabular}{cccccccc}
		\toprule
		& \textbf{Two} & & \textbf{Two} & &  & \textbf{Box} & \textbf{Circle} \\ 
		\textbf{Oval} & \textbf{circles} & \textbf{Box} & \textbf{ovals} & \textbf{Banana} & \textbf{Box2} & \textbf{w/ hole} & \textbf{w/ hole} \\ 
		\midrule
		0.31 & 0.28  & 0.25        & 0.35               & 0.25       & 0.25             & 0.23            & 0.25                \\ \bottomrule
	\end{tabular}%
	\caption{The values of the hyperparameter $\gamma$ for the eight case studies. The hyperparameters are selected based on the incremental approach described in Section~\ref{sec:Lear_validity_domain_Background_One_Class_SVM}.}
	\label{tab:Lear_validity_domain_GAMMA}
\end{table}
The learned boundaries for the case studies are depicted in Figure~\ref{fig:LVD_Bounds_of_illustrative_example}.
As expected, the convex hull does not model holes and disjunct data clusters.
Instead, the convex hull overestimates the validity domain of the case studies.
In contrast, the one-class SVM is able to model holes and disjunct data clusters.
Furthermore, the convex hull includes all data points whereas the one-class SVM also allows for outliers in the data and excludes regions with only small data density from the validity domain (e.g., see the ``box'' case study). 
\subsection{Optimization results}
We minimize the prediction of the eight trained ANNs subject to the the convex hull or the one-class SVM as constraints.
Table~\ref{tab:LVD_optimal_solution_points} shows the optimal solution points, $\textbf{x}^*$, and objective function values, $f_{\text{ANN}}(\textbf{x}^*)$, for the problem with convex hull and one-class SVM as constraints. \newline \indent
\begin{table}[]
	\centering
	\begin{tabular}{lrrrrrrrr}
		\midrule
		\textbf{Case} 	& \multicolumn{3}{c}{\textbf{Convex hull}} & \multicolumn{3}{c}{\textbf{One-class SVM}} & \multicolumn{1}{c}{\textbf{Reference}} \\
		\textbf{study} & \textbf{$x^*$} &  \textbf{$f_{\text{ANN}}^*$} & \textbf{$\Delta$} & \textbf{$x^*$} &  \textbf{$f_{\text{ANN}}^*$} & \textbf{$\Delta$} &  \textbf{$x^*$} \\
		\toprule
		Banana & $(0.1, 2.0)$ & $-5.2$ & $8.52$ & $(0.3, -1.6)$ & $-4.3$ & $0.11$ & $(0.2, -1.6)$\\
		Two circles & $(-1.5, 1.5)$ & $-5.1$ & $3.65$ & $(1.0, 3.7)$ & $-4.8$ & $0.01$ & $(1.3, 3.7)$ \\
		Box & $(0.2, -1.9)$ & $-5.4$ & $2.14$ & $(0.3, -1.5)$ & $-4.5$ &$0.07$ & $(0.2, -1.5)$ \\
		Box w/ hole & $(0.2, 3.6)$ & $-6.4$ & $1.65$ & $(0.5, 3.2)$ & $-4.8$ & $0.70$  & $(0.2, -1.6)$ \\
		Circle w/ hole & $(-1.5, 3.5)$ & $-5.0$ & $0.41$ & $(0.1, 3.6)$ & $-4.6$ & $0.02$ & $(0.2, 3.6)$\\
		Two ovals & $(-1.3, 0.4)$ & $-3.5$ & $0.13$ & $(-1.3, 0.4)$ & $-3.5$ & $0.13$ & $(-1.3, 0.4)$ \\
		Oval & $(1.3, 3.5)$ & $-4.6$  & $0.13$ & $(0.2, -1.6 )$ & $-4.3$ & $0.14$ & $(0.2, -1.6)$ \\
		Box2 & $(0.1, -1.7)$ & $-4.2$ & $0.05$ & $(2.8, 2.9)$ & $-3.8$ &$0.05$ & $(2.9, 2.9)$ \\
		\bottomrule
	\end{tabular}%
	\caption{The table compares the global optimal solutions with convex hull and one-class SVMs (SVMs) as constraints.
		The optimal solution point is given by $x^*$ and the objective function value is given by $f_{\text{ANN}}^*$.
		Also, we provide the error of the data-driven model at the optimal solution $\Delta=\vert f_{\text{ANN}}^*-f_{\text{Peaks}}(x^*)\vert$.
		The reference solution point shows the optimal solution when considering the underlying function $f_{\text{Peaks}}$ as the objective with the one-class SVM constraint.}
	\label{tab:LVD_optimal_solution_points}
\end{table}
The optimal objective function values of the convex hull approach are lower than the ones with the one-class SVM for all case studies because the convex hull overestimates the validity domain.
This overestimation can lead to large errors at the optimal solution points. 
For the ``banana'' case study, the optimal solution found by the convex hull approach is outside the validity domain but at the boundary of the convex hull (see Table~\ref{tab:LVD_optimal_solution_points}).
This leads to a wrongly estimated objective values by the ANN of $-5.2$ with an absolute error of $8.52$.
In contrast, the one-class SVM models the validity domain accurately and yields an optimal solution of $-4.3$ with an absolute error of $0.11$. 
Also, the solution point found by the ANN model with the one-class SVM constraint is close to the reference solution where the learned peaks function is optimized subject to the SVM constraint. 
Similarly, the one-class SVM leads to more reliable results in the data sets ``two circles'', ``box w/ holes'', and ``circle w/ holes''.
Interestingly, the convex hull approach also leads to a substantial prediction error in the ``box'' case study while the one-class SVM models the validity domain accurately. 
This highlights the risk of using the convex hull approach in the presence of noise. \newline \indent 
In Table~\ref{tab:LVD_CPU_Time_OCSVM}, we provide the CPU times for optimization with one-class SVMs embedded.
Using the FS formulation, BARON and MAiNGO perform similarly and solve most problems in the a few hundred CPU seconds. 
The RS formulation outperforms the FS formulation on all problem instances.
In BARON, the speedup factor between the RS and the FS formulation ranges from 5 to over 14.
In comparison, the the speedup factor between the RS and FS in MAiNGO ranges between 583 to over 3,226.
This is in agreement with our previous studies where the McCormick relaxations in the RS lead to smaller subproblems compared to the auxiliary variable method (see Section~\ref{sec:Lear_validity_domain_Method_Optimization}). \newline \indent 
\begin{table}
	\centering
	\begin{tabular}{lrrrrrrr}
		\midrule
		\textbf{Case} &		 \textbf{\# Sup.} 				& \multicolumn{3}{c}{\textbf{BARON}} & \multicolumn{3}{c}{\textbf{MAiNGO}} \\
		\textbf{study} & \textbf{vec.} &  \textbf{FS} & \textbf{RS} &  \textbf{sp-f}  &  \textbf{FS} & \textbf{RS} & \textbf{sp-f} \\
		\toprule
		Oval & 48 & 158.2 s & 35.0 s & \textbf{5} &  197.2 s& 0.20 s& \textbf{986}\\
		Two circles & 50 & 489.5 s& 36.7 s& \textbf{13} & 247.0 s & 0.25 s& \textbf{988}\\
		Box & 52  & 346.6 s& 25.8 s & \textbf{13} & 139.9 s& 0.24 s& \textbf{583} \\
		Two ovals & 58 & 341.3 s & 38.6 s& \textbf{9} & 433.2 s& 0.22 s& \textbf{1,969} \\			
		Banana & 62 & 1,000.0 s& 70.6 s& \textbf{$>$14} & 416.7 s & 0.19 s& \textbf{2,193} \\
		Box2 & 67 & 497.1 s& 22.9 s & \textbf{22} & 1,000.0 s& 0.31 s& \textbf{$>$3,226} \\
		Box w/ hole & 81 & 1,000.0 s & 73.1 s& \textbf{$>$14} & 656.5 s & 0.36 s& \textbf{1,823} \\
		Circle w/ hole & 103 & 1,000.0 s & 75.8 s& \textbf{$>$13} & 1,000.0 s& 0.75 s & \textbf{$>$1,333}\\
		\bottomrule
	\end{tabular}%
	\caption{CPU times for optimization of the eight case studies with the one-class SVM as a constraint.
		The table compares the FS and RS formulations for the BARON and MAiNGO solvers.
		Here, the data-driven model (i.e., the ANN) and the one-class SVM are formulated in the RS and FS.
		The speed up factor (sp-f) is given as the ratio between the FS and the RS solution times.
		Also, the number of support vectors (\# Sup. vec.) is shown as a measure for the problem complexity.}
		\label{tab:LVD_CPU_Time_OCSVM}
\end{table}
In Table~\ref{tab:LVD_CPU_Time_Convex_Hull}, we compare the computational performance for optimization with the convex hull embedded.
The CPU times with the convex hull are lower compared to the ones with one-class SVMs.
MAiNGO is substantially faster than BARON when formulating the problem in the FS.
On average, BARON requires about 83 seconds to solve the problem in the FS while MAiNGO requires only 3 seconds.
The RS formulation again outperforms the FS formulation for all problems.
However, in this case, the speedup factors are in the same order of magnitude for BARON and MAiNGO ranging between 13 and 54.
It should be noted that the RS and the FS formulation of the convex hull constrains are identical.
Therefore, the difference is only due to the formulation of the data-driven model in the objective function, i.e., the ANN (c.f.~\citep{Schweidtmann2019detGlobalANN}). \newline \indent
\begin{table}
	\centering
	\begin{tabular}{lrrrrrrr}
		\midrule
		\textbf{Case} &		 \textbf{\# Sup.} 				& \multicolumn{3}{c}{\textbf{BARON}} & \multicolumn{3}{c}{\textbf{MAiNGO}} \\
		\textbf{study} & \textbf{vec.} &  \textbf{FS} & \textbf{RS} &  \textbf{sp-f}  &  \textbf{FS} & \textbf{RS} & \textbf{sp-f} \\
		\toprule
		Oval & 48 & 74.9 s& 2.8 s& \textbf{27} & 2.7 s & 0.05 s & \textbf{54} \\
		Two circles & 50 & 46.4 s & 3.5 s & \textbf{13} & 2.4 s & 0.09 s& \textbf{27}\\
		Box & 52 & 85.8 s &1.5 s & \textbf{57} & 1.7 s & 0.06 s& \textbf{28} \\
		Two ovals & 58 & 103.9 s&3.1 s& \textbf{34} & 3.8 s & 0.08 s & \textbf{48}\\			
		Banana & 62 & 136.4 s & 5.1 s& \textbf{27} & 3.0 s & 0.09 s& \textbf{33} \\
		Box2 & 67 & 32.3 s & 2.3 s & \textbf{14} & 2.7 s & 0.11 s & \textbf{25}\\
		Box w/ hole & 81 & 51.3 s & 3.5 s & \textbf{15} & 2.4 s & 0.08 s & \textbf{30}  \\
		Circle w/ hole & 103 & 130.2 s & 4.1 s & \textbf{32} & 3.1 s& 0.11 s& \textbf{28}\\
		\bottomrule
	\end{tabular}%
	\caption{CPU times for optimization of the eight case studies with the convex hull as a constraint.
		The table compares the FS and RS formulations for the BARON and MAiNGO solvers. 
		The speed up factor (sp-f) is given as the ratio between the FS and the RS solution times.
		Also, the number of support vectors (\# Sup. vec.) is shown as a measure for the problem complexity.}
		\label{tab:LVD_CPU_Time_Convex_Hull}
\end{table}
\section{Engineering application}
\label{sec:Lear_validity_domain_engineering_case_study}
\begin{figure}
	\centering
	\includegraphics[width=\textwidth]{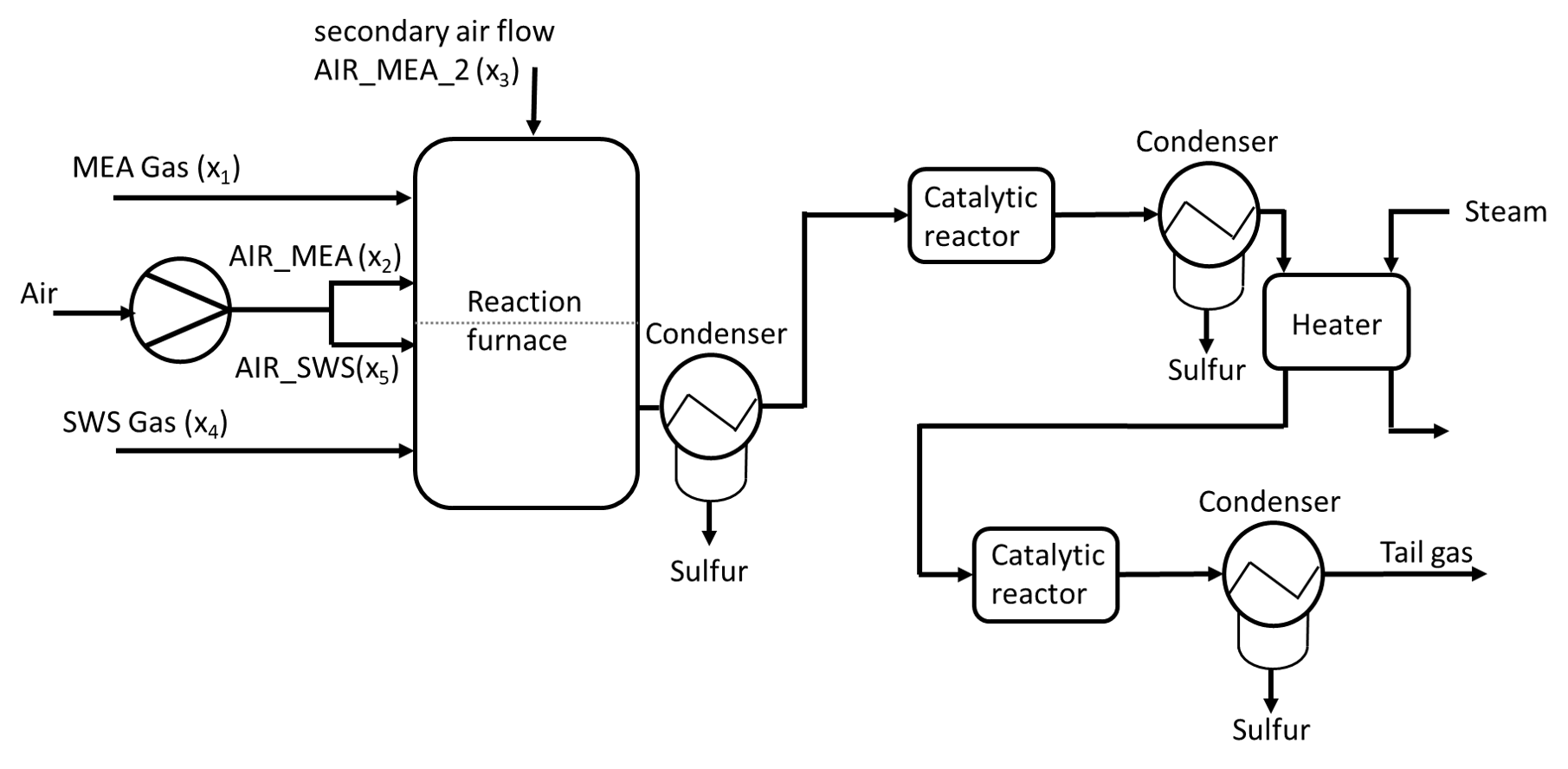}
	\caption{Flowchart of the sulfur recovery unit process}
	\label{fig:LFR_Flowschart_SRU}
\end{figure}

\begin{figure} 
	\centering
	\includegraphics[width=0.5\textheight]{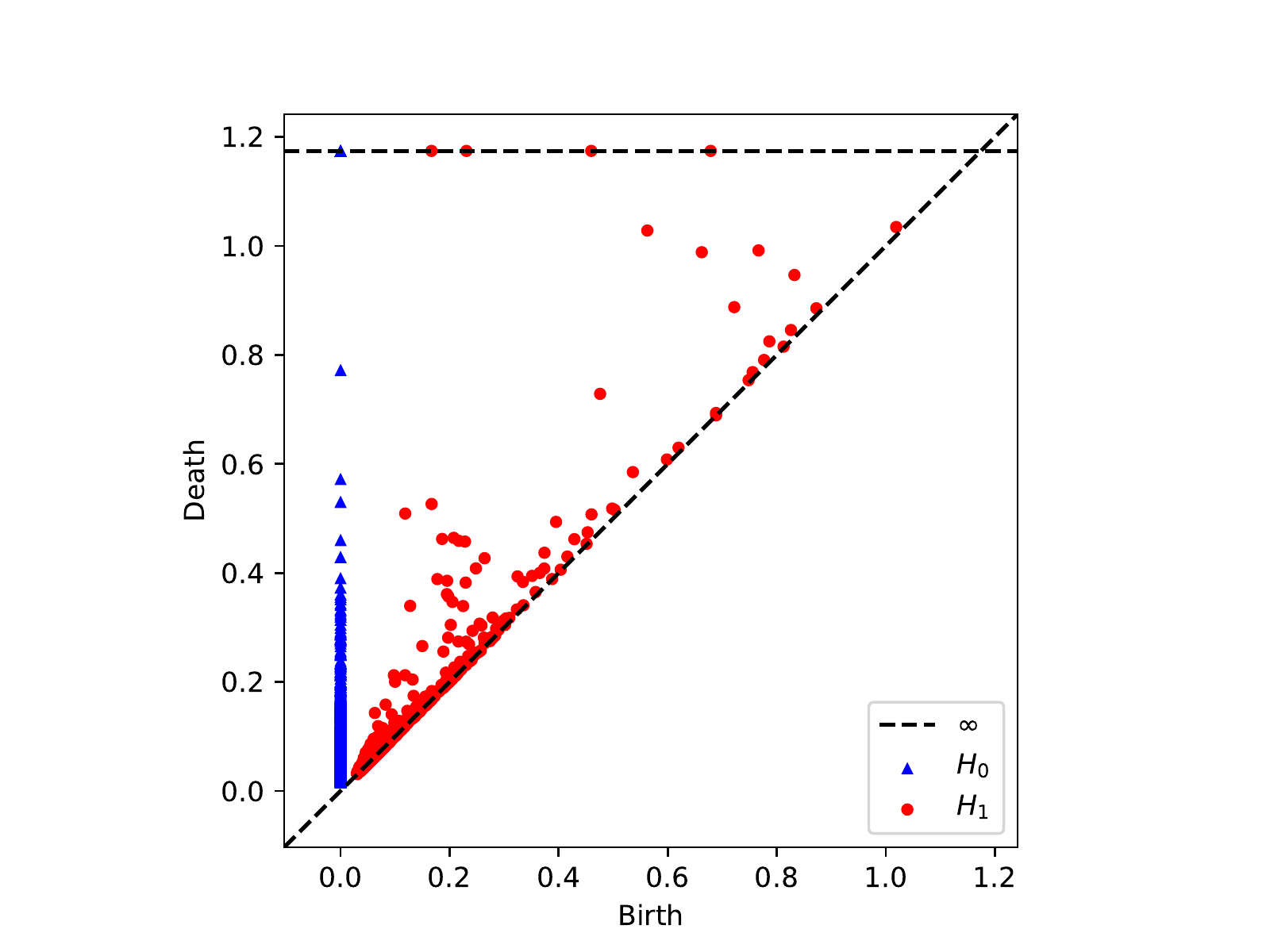}
	\caption{Persistent diagram for of the training data of the engineering case study presented in Section~\ref{sec:Lear_validity_domain_engineering_case_study}}
	\label{fig:Lear_validity_domain_Persistent_diagram_engineering_case_study}
\end{figure}

We consider a sulfur recovery unit as a relevant engineering case study for our work because a large data set of industry operating data is available online for this process.
The efficient recovery of sulfur in petroleum refineries from tail gas is important for environmental reasons.
The sulfur recovery unit process is illustrated in Figure~\ref{fig:LFR_Flowschart_SRU}.
The process has two acid gases as inputs: the MAE gas stream is rich in hydrogen sulphide (H$_2$S) and comes from the gas washing plants.
The SWS gas stream is rich in H$_2$S and ammonia and comes from a sour water stripping plant.
In the sulfur recovery unit, the acid gases are burnt via partial reaction with air in a two-chamber reaction furnace.
Then, the combustion product is further treated in two subsequent catalytic reactors resulting in a tail gas stream that contains residuals of H$_2$S and sulfur dioxide (SO$_2$). 
A detailed process description can be found in the literature~\citep{fortuna2007soft}. \newline \indent
A key issue of the sulfur recovery unit is the control of the secondary air flow to ensure optimal conditions for the total removal of the sulfur compounds in the catalytic converters.
Previous works have investigated soft sensors for the tail gas concentrations of hydrogen sulphide (H$_2$S) and sulfur dioxide (SO$_2$) using ANNs and implemented those in industry for monitoring~\citep{quek2000consider,fortuna2003soft,fortuna2007soft}.
In this case study, we solve an open-loop control problem to find the optimal secondary air flow rate.
The objective is to minimize $\vert c_{H_2S}-2\cdot c_{SO_2} \vert$ such that the two reactants are in stoichiometric proportion.
Similar to the previous literature by \citet{fortuna2003soft,fortuna2007soft}, we also train two ANNs to predict the concentrations:
\begin{align*}
&c_{H_2S, k} &= f_{{ANN, H_2S}} \left( x_{1,k}, x_{1,k-5}, x_{1,k-7}, x_{1,k-9}, ... x_{5,k}, x_{5,k-5}, x_{5,k-7}, x_{5,k-9} \right), \\
&c_{SO_2, k} &= f_{{ANN, SO_2}} \left( x_{1,k}, x_{1,k-5}, x_{1,k-7}, x_{1,k-9}, ... x_{5,k}, x_{5,k-5}, x_{5,k-7}, x_{5,k-9} \right),
\end{align*}
where $x_{1,k}$ is the gas flow in the MEA zone, 
$x_{2,k}$ is the air flow in the MEA zone, 
$x_{3,k}$ is the secondary air flow in the MEA zone, 
$x_{4,k}$ is the air flow in the SWS zone, 
$x_{5,k}$ is the gas flow in the SWS zone
at time step $k$.
The $SO_2$ ANN has one hidden layer with eight neurons and the $H_2S$ ANN has two hidden layers with eight neurons each.
The data-driven models are trained on (scaled) industrial data collected at a plant located in Priolo, Italy available at \url{https://www.openml.org/d/23515}. 
The data set includes a time series with approximately 10,000 data samples and we use the first 90\% of the data for training.
The control variable of the NMPC is the secondary air flow $x_{3,k}$ while the other inputs are observable parameters.
As the control is critical for process safety, the validity limits of the data-driven model should be considered. \newline \indent
In order to analyze the topology of the 20-dimensional input training data set of ANNs, we perform persistent homology.
Due to the large number of data points, the exact computation of the persistent diagram is expensive.
We apply approximate sparse filtration instead~\citep{cavanna2015geometric}.
The persistent diagram for this case study is shown in Figure~\ref{fig:Lear_validity_domain_Persistent_diagram_engineering_case_study}.
The diagram shows that there exist a number of holes in the data set that persist over a long time span.
Also, a separate cluster can be observed in the data.
This motivates the use of one-class SVM to obey validity limits of data-driven models. \newline \indent
As we have no physical model of the process available, the closed-loop performance of the controller is not studied in this example.
For illustration, we perform one step of an open-loop controller for the secondary air flow $x_{3,k}$. 
We select a random operating point from the historic plant data (Table~\ref{tab:LVD_SRU}) and let the solver determine the optimal control action.
The problem is solved to global optimality within 0.33 CPU seconds and identifies a control action $x_{3,k}=0.266$ that results in the desired stoichiometric composition, i.e., $\vert c_{H_2S}-2\cdot c_{SO_2} \vert=1.3\cdot 10^{-5}$.
This engineering case study also demonstrates the potential of the proposed method for NMPC. 
Note that deterministic global NMPC can become computationally expensive for long control horizons and higher dimensional control vectors~\citep{chachuat2006global,Doncevic2020RNN_Global_Control,mitsos_20_hw}.
\begin{table}[]
	\centering
	\begin{tabular}{lrrrr}
		\midrule
		& $k$ & $k-5$ & $k-7$ & $k-9$ \\
		\toprule
		$x_{1}$ & 0.627 & 0.6215 & 0.623 & 0.622 \\
		$x_{2}$ & 0.770 & 0.769 & 0.754 & 0.769 \\
		$x_{3}$ & $x_{3,k}$ & 0.174 & 0.192 & 0.198 \\
		$x_{4}$ & 0.376 & 0.399 & 0.415 & 0.410 \\
		$x_{5}$ & 0.513 & 0.512 & 0.511 & 0.504 \\
		\bottomrule
	\end{tabular}%
	\caption{Operating point of the sulfur recovery unit that is considered for the NMPC optimization step.}
	\label{tab:LVD_SRU}
\end{table}
\section{Conclusion}
Safety concerns and extrapolation issues often impede industrial applications of machine learning models. 
We present a three-step approach to obey the validity limits of data-driven models. 
First, we perform a data topology analysis using persistent homology. 
Second, we model the validity domain of the data-driven model using either the convex hull or a one-class SVM.
Third, we perform deterministic global optimization with the validity domain model as a constraint. \newline \indent 
All used and developed methods are available open-source.
Also, we currently develop a Python interface for our solver MAiNGO.  
Thus, all methods can be applied and further developed in academia and industry for free. \newline \indent 
Our method has the potential to enhance safety, trust, and reliability of machine learning approaches. 
Moreover, we demonstrate that persistent homology is a valuable method for understanding the topology of data in high dimensional spaces. 
Besides industry applications, promising future work also includes the application to optimization problems occurring in molecular design where molecules are parameterized through graph neural networks~\citep{schweidtmann2020graph} or autoencoders~\citep{jin2018junction}.
Also, time-dependent design space descriptions are desired in pharmaceutics~\citep{von2020working}.
The proposed method can also be extended by considering and comparing other one-class classification methods. 
\section*{Acknowledgements}
	We are grateful to Beno{\^{i}}t Chachuat for providing MC\texttt{++} under Eclipse Public License. 
	We also thank Dominik Bongartz and Jaromi{\l} Najman for their work on MAiNGO. 

\section*{Declarations}
\subsection*{Code availability}
The method is ready-to-use and available open-source as part of our  ``MeLOn - \textbf{M}achin\textbf{e} \textbf{L}earning Models for \textbf{O}ptimizatio\textbf{n}'' toolbox under the Eclipse public license (\url{https://git.rwth-aachen.de/avt.svt/public/MeLOn}).

\subsection*{Availability of data and material}
The (scaled) industrial data used in the engineering case study are available at \url{https://www.openml.org/d/23515}. 

\subsection*{Authors' contributions}
AMS and JMW designed the research concept.
AMS wrote the manuscript.
JMW run the persistent homology analyzed the persistent plots, and wrote the corresponding method and result sections.
AMS, CW, and LN run the optimization.
AMS, LN, and CW implemented the model in the MeLOn tool.
AM is principal investigator who guided the effort and edited the manuscript. 
\subsection*{Conflicts of interest}
The authors declare that they have no conflict of interest.
\subsection*{Funding}
This work was supported by the Deutsche Forschungsgemeinschaft (DFG, German Research Foundation) under Germany's Excellence Strategy - Cluster of Excellence 2186 ``The Fuel Science Center''.
\bibliographystyle{spbasic}      

\bibliography{Bibs}   

\end{document}